\documentclass[preprint]{imsart}
\RequirePackage[OT1]{fontenc}
\RequirePackage{amsthm,amsmath}
\RequirePackage[numbers]{natbib}
\RequirePackage{hypernat}
\usepackage{amssymb,graphicx,url,bbm}
\usepackage[colorlinks,citecolor=blue,urlcolor=blue]{hyperref}

\startlocaldefs
\def\d{\delta}

\def\th{\theta}

\def\P{{\mathbb P}}
\def\G{{\mathbb G}}

\def\a{\alpha}
\def\b{\beta}
\def\R{\mathbb R}
\def\dd{\Delta}

\def\d{\delta}

\def\G{{\mathbb G}}

\def\IK{I\!\!K}

\def\P{{\mathbb P}}

\def\labda1{\lambda_1}
\def\labda2{\lambda_2}

\def\f{\phi}

\def\k{\kappa}

\def\s{\sigma}

\def\comment#1{\relax}

\def\=in{\mathop{\rm =}}

\def\eop{\hfill\mbox{$\Box$}\newline}

\newtheorem{theorem}{Theorem}[section]

\newtheorem{remark}{Remark}[section]

\numberwithin{equation}{section}
\theoremstyle{plain}
\endlocaldefs

\usepackage{amsmath,amssymb,here,amsfonts,latexsym,sym1,graphicx,here}
\usepackage{amsthm,curves}
\usepackage{natbib}
\usepackage{graphicx}
\usepackage{caption}
\usepackage{subcaption}
\usepackage{float}

\advance \topmargin by -\headheight \advance \topmargin by
-\headsep

\def\th{\theta}

\def\P{{\mathbb P}}

\def\G{{\mathbb G}}

\include{commandsCUP}

\begin{document}

\begin{frontmatter}
\title{The bivariate current status model}
\runtitle{Bivariate current status}

\begin{aug}
\author{\fnms{Piet} \snm{Groeneboom}\corref{}\ead[label=e1]{P.Groeneboom@tudelft.nl}
\ead[label=u1,url]{http://dutiosc.twi.tudelft.nl/\textasciitilde pietg/}}
\address{Delft Institute of Applied Mathematics,
Delft University of Technology,\\
Mekelweg 4, 2628 CD Delft,
The Netherlands\\
\printead{e1},
\printead{u1}
}
\affiliation{Delft University}
\runauthor{P.\ Groeneboom}
\end{aug}

\begin{abstract}
For the univariate current status and, more generally, the interval censoring model, distribution theory has been developed for the maximum likelihood estimator (MLE) and smoothed maximum likelihood estimator (SMLE) of the unknown distribution function, see, e.g., \cite{GrWe:92}, \cite{piet:96}, \cite{GeGr:96}, \cite{GeGr:97}, \cite{GeGr:99}, \cite{piet_geurt_birgit:10}, \cite{piet_tom:11} and  \cite{piet:11e}. For the bivariate current status and interval censoring models distribution theory of this type is still absent and even the rate at which we can expect reasonable estimators to converge is unknown.

We define a purely discrete plug-in estimator of the distribution function which locally converges at rate $n^{1/3}$, and derive its (normal) limit distribution. Unlike the MLE or SMLE, this estimator is not a proper distribution function. Since the estimator is purely discrete, it demonstrates that the $n^{1/3}$
convergence rate is in principle possible for the MLE, but whether this actually holds for the MLE is still an open problem. 

We compare the behavior of the plug-in estimator with the behavior of the MLE on a sieve and the SMLE in a simulation study. This indicates that the plug-in estimator and the SMLE have a smaller variance but a larger bias than the sieved MLE. The SMLE is conjectured to have a $n^{1/3}$-rate of convergence if we use bandwidths of order $n^{-1/6}$. We derive its (normal) limit distribution, using this assumption. Finally, we demonstrate the behavior of the MLE and SMLE for the bivariate interval censored data of \cite{Betensky:99}, which have been discussed by many authors, see  e.g., \cite{sun:06}, \cite{gentleman_vandal:02}, \cite{bogaerts:04} and \cite{marloes:2005}.
\end{abstract}



\end{frontmatter}

\section{Introduction}
\label{section:intro}
\setcounter{equation}{0}
We consider the bivariate current status model, also called the bivariate interval censoring, case 1, model. This means that our observations consist of a quadruple $(T,U,\dd_1,\dd_2)$, where
\begin{equation}
\label{delta_bivar}
\dd_1=1_{\{X\le T\}},\,\dd_2=1_{\{Y\le U\}},
\end{equation}
and $(X,Y)$ is independent of the observation $(T,U)$. We want to estimate the distribution function $F_0$ of the `hidden' random vector $(X,Y)$.

A maximum likelihood estimator $\hat F_n$ of $F_0$, the distribution function of $(X,Y)$, maximizes the expression
\begin{align*}
&\sum_{i=1}^n\left\{\dd_{i1}\dd_{i2}\log F(T_i,U_i)
+\dd_{i1}\left(1-\dd_{i2}\right)\log\left\{F(T_i,\infty)-F(T_i,U_i)\right\}
\right.\\
&\qquad+\left(1-\dd_{i1}\right)\dd_{i2}\log\left\{F(\infty,U_i)-F(T_i,U_i)\right\}\\
&\qquad\left.+\left(1-\dd_{i1}\right)\left(1-\dd_{i2}\right)
\log\left\{1-F(\infty,U_i)-F(T_i,\infty)+F(T_i,U_i)\right\}\right\}
\end{align*}
over all bivariate distribution functions $F$. Another formulation is that $\hat F_n$ maximizees
\begin{align*}
&\int \d_1\d_2\log F(u,v)\,d\,\P_n
+\int
\d_1(1-\d_2)\log\left\{F_1(u)-F(u,v)\right\}\,d\,\P_n\\
&\qquad+\int (1-\d_1)\d_2\log\left\{F_2(v)-F(u,v)\right\}
\,d\,\P_n\\
&\qquad+\int (1-\d_1)(1-\d_2)\log\left\{1-F_1(u)-F_2(v)+F(u,v)\right\}\,d\P_n
\end{align*}
over $F$, where $F_1$ and $F_2$ are the first and second marginal dfs of $F$, respectively, and $\P_n$ is the empirical measure of the observations $(T_i,U_i,\dd_{i1},\dd_{i2})$, $i=1,\dots,n$.

One looks for a solution of the form
$$
\hat F_n=\sum_{j=1}^m\a_j1_{[\tau_j,\mathbf\infty)},\,\sum_{j=1}^m\a_j\le1,\,\a_j>0,\,1\le j\le m=m_n,
$$
where we denote by $[\tau_j,\boldmath\infty)$ an infinite rectangle $[t_{j1},\infty)\times[t_{j2},\infty)$, where
$\tau_j=(t_{j1},t_{j2})$. Then the (Fenchel or Kuhn-Tucker) duality conditions for the solution are:
\begin{align}
\label{fenchel-ineq-ML}
&\int_{[{\mathbf t},\infty)} \frac{\d_1\d_2}{F(u,v)}\,d\,\P_n
+\int_{[t_1,\infty)\times[0,t_2)}
\frac{\d_1(1-\d_2)}{F_1(u)-F(u,v)}\,d\,\P_n\\
&+\int_{[0,t_1)\times[t_2,\infty)}\frac{(1-\d_1)\d_2}{F_2(v)-F(u,v)}
\,d\,\P_n+\int_{[0,t_1)\times[0,t_2)}\,
\frac{(1-\d_1)(1-\d_2)}{1-F_1(u)-F_2(v)+F(u,v)}\,d\,\P_n\nonumber\\
&\le 1,\nonumber
\end{align}
for all ${\mathbf t}=(t_1,t_2)\in\R^2$, where $\P_n$ is the empirical df of the observations
$(u_i,v_i,\d_{i1},\d_{i2})$.
We must have equality in (\ref{fenchel-ineq-ML}), if ${\mathbf t}=\tau_j$, $j=1,\dots,m_n$, is a point of mass of the solution (i.e., the constraints are active!), where the rectangles $[\tau_j,\infty)$ are the generators of the solution.

An R package, called `MLEcens' is available for computing the MLE. The algorithm determines the maximal intersection rectangles where the MLE has mass via a preliminary reduction algorithm, and next computes the mass of the MLE in these rectangles, using the support reduction algorithm of \cite{piet_geurt_jon:08}. The reduction algorithm is described in \cite{marloes:2005}. The R package uses as an example a data set, studied in \cite{Betensky:99}, which is actually not of the current status type but has interval censoring, case 2, data. We shall discuss these data in section \ref{section:interval_cens}. The MLE for this data set is also discussed in section 7.3.3 of \cite{sun:06}, who also refers to \cite{gentleman_vandal:02} and \cite{bogaerts:04} for discussions of the computation of the MLE for this data set. The computation of the rectangles where the MLE puts mass is also treated in \cite{song:01}, where also  minimax lower bounds for the estimation rate of the MLE and consistency of the MLE are derived.

There is an extensive discussion on where to put the mass, once one has determined rectangles which can have positive mass, see, e.g. \cite{sun:06}, section 7.3, \cite{gentleman_vandal:02}, \cite{bogaerts:04} and \cite{marloes:2005}. The algorithm for computing these rectangles, proposed in \cite{marloes:2005}, seems at present to be the fastest.

It is in our view somewhat doubtful whether all the energy spent on computing these maximal intersection rectangles and the ensuing question of whether one should place the mass of the MLE at the left lower corner or the right upper corner of the rectangles is really worth the effort. One could also specify in advance a set of points where one allows mass to be placed. In this way one obtains an MLE on a sieve, where the sieve consists of distributions having discrete mass at these points. The bottleneck in the computation of the MLE for the bivariate interval censoring problem is not the determination of the maximal intersection rectangles, but the computation of the mass the MLE puts on these rectangles, since there usually are very many!

The latter phenomenon shows up in particular in simulations. As an example, simulating data from the distribution with density $f(x,y)=x+y$ on the unit square, with a uniform observation distribution, we got for sample size $n=5000$ about $5\cdot 10^5$ possible rectangles where the masses could be placed, which is (at present) an almost prohibitive number if one wants to do simulations of the limit behavior of the MLE on an ordinary table computer or laptop.
Ultimately, the discussion on these matters should in our view be determined by insights into the distribution theory of the MLE or the MLE on the chosen sieve.

In section \ref{sec:asymp_bivarCS} we show that a purely discrete plug-in estimator locally attains the $n^{1/3}$ rate. We also determine its asymptotic (normal) distribution. In section \ref{sec:SMLE} we study the local limit behavior of the smoothed maximum likelihood estimator and derive its asymptotic distribution under the assumption that it can asymptotically can be represented by an integral in the observation space, proceeding along similar lines as in the one-dimensional case (see, e.g., \cite{GeGr:97}, \cite{piet_geurt_birgit:10} and \cite{piet:12b}). Section \ref{section:simulations} presents a simulation study, comparing the behavior of the MLE, the SMLE and the plug-in estimator. The results seem to be in accordance with Theorems \ref{th:plug-in} and \ref{th:SMLE} in sections \ref{sec:asymp_bivarCS} and \ref{sec:SMLE}, respectively. Section \ref{section:interval_cens} extends the MLE and SMLE to a more general setting of interval censoring and applies the methods on a data set in \cite{Betensky:99}, which has been discussed by many authors, see  e.g., \cite{sun:06}, \cite{gentleman_vandal:02}, \cite{bogaerts:04} and \cite{marloes:2005}. The paper ends with some concluding remarks on faster rates for the SMLE, which can be attained if one uses higher order kernels.

\section{A purely discrete $n^{1/3}$-rate estimator for the bivariate current status model}
\label{sec:asymp_bivarCS}
Basically, the MLE for the $1$-dimensional current status model is the monotone derivative of the cusum diagram
$$
\left(\G_n(t),V_n(t)\right),\,t\in I,\,\qquad V_n(t)=\int_{u\le t}\d\,d\P_n(u,\d),
$$
where $I$ is the observation interval and $\G_n$ the empirical distribution function of the observations $T_1,\dots,T_n$. So it can be considered to be a monotone version of the `derivative' $dV_n(t)/d\G_n(t)$. Note that if we replace $V_n$ and $\G_n$ by their deterministic equivalents, the derivative becomes
$$
\frac{F_0(t)g(t)}{g(t)}=F_0(t),
$$
so is indeed the object we want to estimate.

For the simplest bivariate current status model, which is sometimes called the `in-out' model, we only have the information whether the hidden variable is below and to the left of the observation point $(T_i,U_i)$ or not. In this case we could also define
$$
V_n(t,u)=\int_{v\le t,\,w\le u}\d\,d\P_n(v,w,\d),
$$
where $\d=1$ represents the situation that the hidden variable is below and to the left of $(v,w)$. If the empirical observation distribution is again denoted by $\G_n$, we this time want to estimate the `derivative' $dV_n(t,u)/d\G_n(t,u)$, since, replacing $V_n$ and $\G_n$ by their deterministic equivalents, the derivative becomes
$$
\frac{F_0(t,u)g(t,u)}{g(t,u)}=F_0(t,u).
$$
So we want to find a version of the derivative $dV_n(t,u)/d\G_n(t,u)$, under the (shape) restriction that it is a bivariate distribution function.

However, a natural cusum diagram for this situation does not seem to exist. But we can define a $2$-dimensional `Fenchel process', incorporating the duality conditions for a solution of the optimization problem. Analogously to the $1$-dimensional current status model, the Fenchel duality conditions for the isotonic least squares (LS)estimate, minimizing
$$
\sum_{i=1}^n\left\{\dd_i-F(T_i,U_i)\right\}^2,\qquad \dd_i=1_{\{X_i\le T_i,\,Y_i\le U_i\}}
$$
over all bivariate distribution functions $F$, where the $(X_i,Y_i)$ are the hidden variables, are:
\begin{equation}
\label{fenchel_in_out}
\int_{v\ge t,\,w\ge u}\left\{\d-F(v,w)\right\}\,d\P_n(v,w,\d)\le0,
\end{equation}
with equality if $(t,u)$ is a point of mass of the solution.  So we have to deal with a process
\begin{equation}
\label{invelope_simpleCS_bivar}
(t,u)\mapsto \int_{v\ge t,\,w\ge u}F(v,w)\,d\G_n(v,w)
\end{equation}
which has to lie above the process
$$
(t,u)\mapsto \int_{v\ge t,\,w\ge u}\d\,d\P_n(v,w,\d),
$$
with points of touch at points of mass of $F$. Denoting temporarily the process (\ref{invelope_simpleCS_bivar})
by $Q_n$, we get that the isotonic least squares estimator (which no longer necessarily coincides with the MLE!) can (formally) be denoted by $dQ_n(t,u)/d\G_n(t,u)$. Note, however, that the function $Q_n$ is not necessarily close to a convex or concave function, so here the analogy with $1$-dimensional current status breaks down. But it must have the property that its `derivative' w.r.t.\ $d\G_n$ must be a distribution function, which is analogous to the fact that the derivative of the convex minorant of the cusum diagram must be a distribution function in the $1$-dimensional case.

For the real current status model the situation is more complicated, since we then have to deal with $4$ regions instead of $2$. From (\ref{fenchel-ineq-ML}) we get:
\begin{align*}
&\int_{[{\mathbf t},\infty)} \frac{\d_1\d_2}{\hat F_n(u,v)}\,d\,\P_n
+\int_{[t_1,\infty)\times[0,t_2)}
\frac{\d_1(1-\d_2)}{\hat F_{n1}(u)-\hat F_n(u,v)}\,d\,\P_n\nonumber\\
&+\int_{[0,t_1)\times[t_2,\infty)}\frac{(1-\d_1)\d_2}{\hat F_{n2}(v)-\hat F_n(u,v)}
\,d\,\P_n+\int_{[0,t_1)\times[0,t_2)}\,
\frac{(1-\d_1)(1-\d_2)}{1-\hat F_{n1}(u)-\hat F_{n2}(v)+\hat F_n(u,v)}\,d\,\P_n\nonumber\\
&\le 1,
\end{align*}
where ${\mathbf t}=(t_1,t_2)$, with equality if $(t_1,t_2)$ is a point of mass of $\hat F_n$. 

It has been conjectured that the MLE in the bivariate current status model converges locally at rate $n^{1/3}$, just as in the $1$-dimensional current status model (with smooth underlying distribution functions). \cite{song:01} proves a minimax lower bound of order $n^{-1/3}$. It would be somewhat surprising if the $1$-dimensional rate would be preserved in dimension $2$, since in general one gets lower rates for density estimators if the dimension gets up, and the estimation of the distribution function in the current status model is similar to density estimation problems, as argued above. 

To show that it is in principle possible to attain the local rate $n^{1/3}$, we construct a purely discrete estimator, converging locally at rate $n^{1/3}$. We restrict ourselves for simplicity to distributions with support $[0,1]^2$ in the remainder of this section, but the generalization to more general rectangles is obvious. We have the following result, which is proved in the Appendix.

\begin{theorem}
\label{th:plug-in}
Consider an interior point $(t,u)$, and define the square $A_n$, with midpoint $(t,u)$, by:
$$
A_n=[t-h_n,t+h_n]\times[u-h_n,u+h_n].
$$
Moreover, suppose that the observation distribution $G$ is twice continuously differentiable at $(t,u)$ with a strictly positive density $g(t,u)$ at $(t,u)$, and that $F_0$ is twice continuously differentiable at $(t,u)$. Moreover, suppose
\begin{equation}
\lim_{n\to\infty}h_n^2n^{1/3}=c>0.
\end{equation}

Then the estimator
\begin{align}
\label{plug-in_bivarCS}
\tilde F_n(t,u)
     & \stackrel{\mbox{\small def}} = \frac{\int_{A_n} \d_1\d_2\,d\P_n(v,w,\d_1,\d_2)}{\int_{A_n} d \G_n(v,w)}\,, 
\end{align}
where $\G_n$ is the empirical distribution function of the observations $(T_i,U_i)$ and $\P_n$ is the empirical distribution function of the observations
$$
\left(T_i,U_i,\dd_{i1},\dd_{i2}\right),\,i=1,\dots,n,
$$
satisfies:
$$
n^{1/3}\left\{\tilde F_n(t,u)-F_0(t,u)\right\}\stackrel{{\cal D}}\longrightarrow
N\left(\beta,\s^2\right),
$$
where $N(\b,\s^2)$ is a normal distribution with first moment
$$
\b=c\left\{\tfrac16\left\{\partial_1^2F_0(t,u)+\partial_2^2F_0(t,u)\right\}
+\frac{\partial_1F_0(t,u)\partial_1g(t,u)+\partial_2F_0(t,u)\partial_2g(t,u)}{3g(t,u)}\right\},
$$
and variance
\begin{align*}
\s^2=\frac{F_0(t,u)\left\{1-F_0(t,u)\right\}}{4cg(t,u)}\,.
\end{align*}
\end{theorem}

We now allow as possible points of mass the points $(t_{n1},u_{n1}),\dots, (t_{n,m_n},u_{n,m_n})$ running through a rectangular grid, where the distances between the points on the $x$- and $y$-axis are of order $n^{-1/3}$, and define the estimate $\tilde F_n$ at each point $(t_{ni},u_{ni})$ as in Theorem \ref{th:plug-in}. Next we define the masses $p_{ni}$ at the points $(t_{ni},u_{ni})$ by the equations
$$
\sum_{j:t_{nj}\le t_{ni},\,u_{nj}\le u_{ni}}p_{nj}=\tilde F_n(t_{ni},u_{ni}),\,i=1,\dots,m_n.
$$
Note that the estimate $\tilde F_n$ we obtain in this way is not necessarily a distribution function and that the masses $p_{nj}$ can have negative values.

\begin{figure}[!ht]
\begin{subfigure}[b]{0.45\textwidth}
\centering
\includegraphics[width=\textwidth]{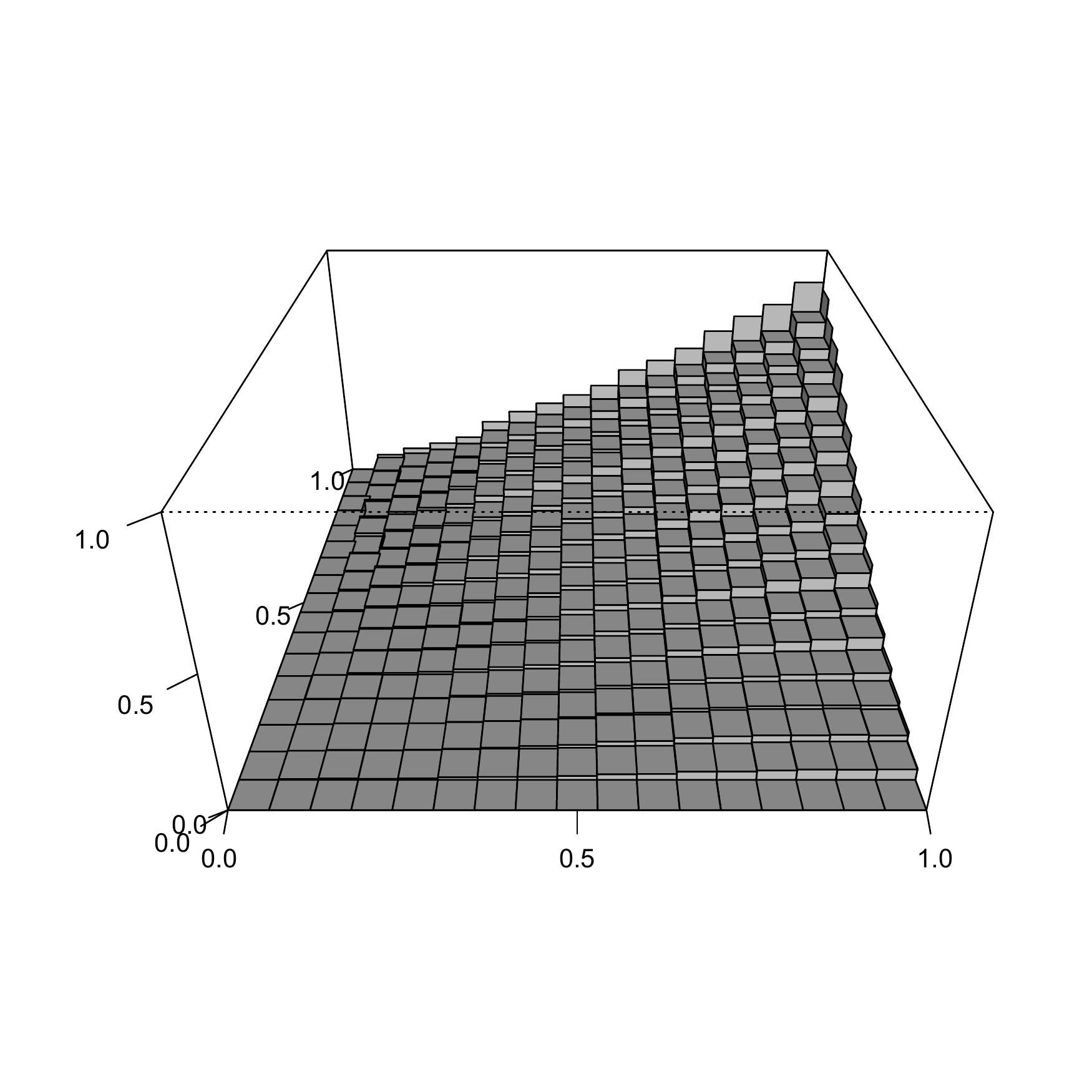}
\caption{Plug-in estimator}
\label{fig:plug-in_bivar}
\end{subfigure}
\begin{subfigure}[b]{0.45\textwidth}
\centering
\includegraphics[width=\textwidth]{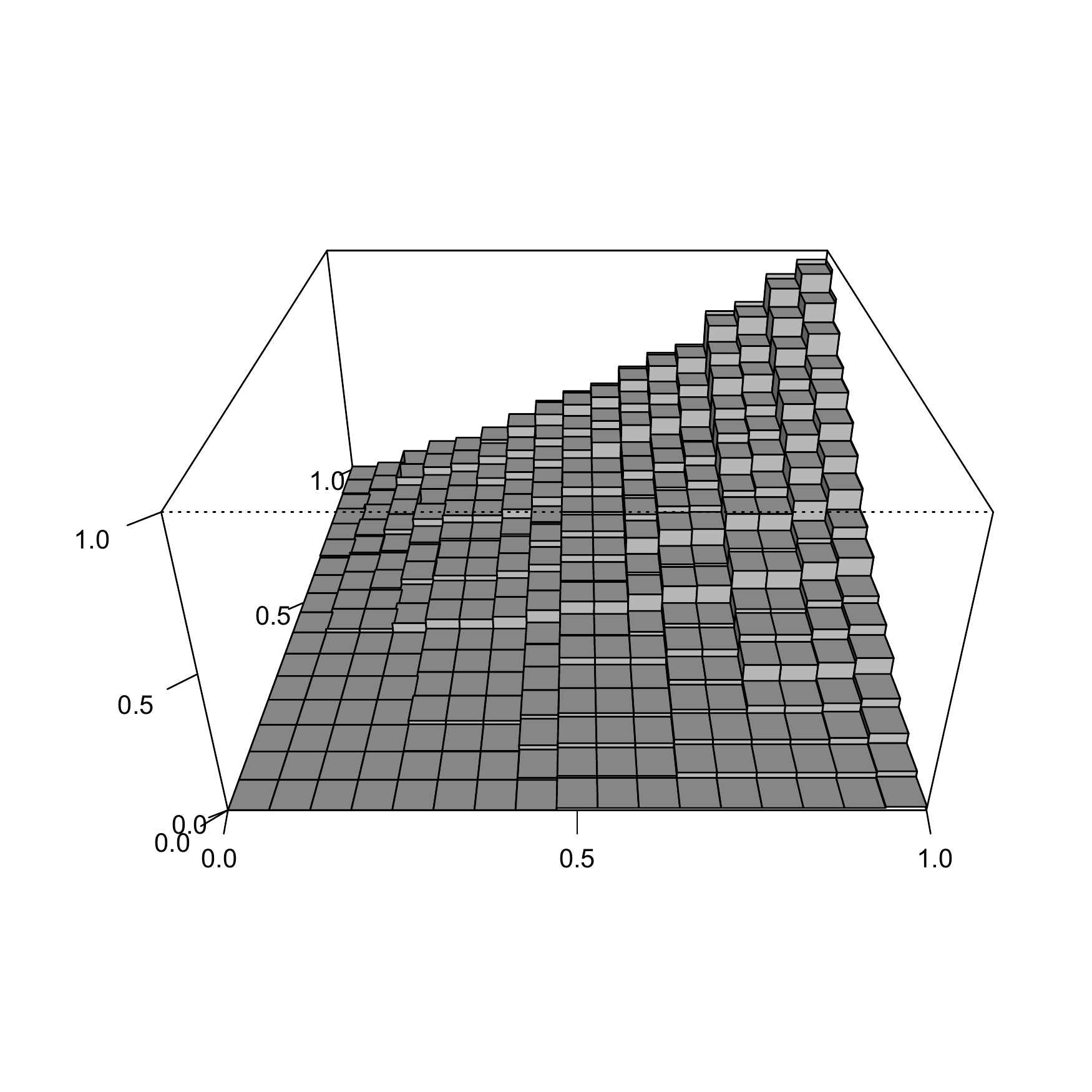}
\caption{MLE on the points of mass of the plug-in estimate}
\label{fig:MLE_bivar}
\end{subfigure}
\caption{The plug-in estimator $\tilde F_n$ and the MLE, for a sample of size $n=5000$ of bivariate current status data, where the hidden variables have a distribution with density $f_0(x,y)=x+y$, and the observation distribution is uniform on $[0,1]^2$.}
\label{fig:plug-in+MLE_bivarCS}
\end{figure}

Also note that we get roughly order $n^{1/3}\times n^{1/3}$ equations in this way, which turned out to be solvable, although it is not clear beforehand that the system one gets is non-singular. Nevertheless, one can build up the system from left and below up to the right and above, where one gets more an more values in the corresponding matrix, so it seems likely that in general the solution exists. This seems a point for further research.
A picture of $\tilde F_n$, together with the MLE, computed on the sieve of points of mass of the plug-in estimator, is shown in Figure \ref{fig:plug-in+MLE_bivarCS}. The sieved MLE is a proper (discrete) distribution function, so all its masses are nonnegative.

Assuming the support of the distribution with function $F_0$ to be $[0,1]^2$, we treat the points near the upper and right boundary in a special way in computing $\tilde F_n$. If, for example $t_{ni}>1-h_n$, where $h_n\sim n^{-1/6}$ and $(t_{ni},u_{ni})$ is a point of the grid, we put the $\d_{j1}$ corresponding to the contribution of observations $(T_j,U_j)$ such that $T_j\ge 2-t_{ni}-h_n$, equal to $1$. We treat the second coordinate in the same way. This reduces the bias we otherwise would get, with an underestimation of the distribution function near the right and upper boundary. The idea to treat the points near the boundary in this way is inspired by, but different from, the reflection method proposed by \cite{schuster:85}. For points near the left and lower boundary, we follow a similar procedure. If, for example $t<h_n$, where $h_n\sim n^{-1/6}$ and $(t,u)$ is a point of the grid, we put the $\d_{j1}$ corresponding to the contribution of observations $(T_j,U_j)$ such that $T_j\le h_n-t$, equal to $0$. The bias near the boundary is actually of order $O(h_n)$ in this way and we do not attain the $O(h_n^2)$ for interior points, though.

\section{The smoothed maximum likelihood estimator}
\label{sec:SMLE}
Throughout this section we will assume for simplicity that the support of the distribution of the `unobservables' is $[0,1]^2$ and that we want to estimate the corresponding distribution function $F_0$ on $[0,1]^2$. The generalization to more general rectangles is obvious.

Let $K$ be a symmetric non-negative kernel, for example the triweight kernel
$$
K(x)=\tfrac{35}{32}\left(1-x^2\right)^31_{[-1,1]}(x),
$$
and $K_h(x)=h^{-1}K(x/h)$, for $h>0$. Moreover, let the integrated kernel $\IK$ be defined by
$$
\IK(x)=\int_{-\infty}^x K(y)\,dy.
$$
We follow the approach for the $1$-dimensional case, discussed in the references \cite{GeGr:96} to \cite{piet_geurt_birgit:10}.

At an interior point $(t,u)$, not too close to the boundary, the smoothed maximum likelihood estimator (SMLE) is just defined by
\begin{equation}
\label{SMLE}
\hat F_{nh}^{(SML)}(t,u)=\int \IK_h(t-v)\IK_h(u-w)\,d\hat F_n(v,w),\qquad \IK_h(x)=\IK(x/h).
\end{equation}
To prevent the negative bias at the right and upper boundary of the support, we also perform a correction by  extending the definition of $\IK$ near the upper and right boundary by:
\begin{equation}
\label{boundary_IK}
\IK^b(x)=\int_{x}^{\infty} K(y)\,dy
\end{equation}
and defining
\begin{equation}
\label{SMLE2}
\hat F_{nh}^{(SML)}(t,u)=\int\left\{\IK_h(t-v)+\IK_h^b(2-t-v)\right\}
\left\{\IK_h(u-w)+\IK_h^b(2-u-w)\right\}\,d\hat F_n(v,w),
\end{equation}
and $\IK_h^b(x)=h^{-1}\IK^b(x/h)$. This definition of the (integrated) boundary kernel is based on the reflection boundary correction method for density estimates, proposed in \cite{schuster:85}. Note that the definitions (\ref{SMLE}) and (\ref{SMLE2}) coincide if $\max(t,u)\le 1-h$.

We next define the score function in the hidden space:
$$
\k_{(t,u)}(x,y)=\left\{\IK_h(t-x)+\IK_h^b(2-t-y)\right\}
\left\{\IK_h(u-x)+\IK_h^b(2-u-y)\right\},
$$
Scores in the observation space are given by
\begin{equation}
\label{score_definition}
\th_{F_0}(v,w,\d_1,\d_2)=E\left\{a(X,Y)\bigm|\left(T,U,\dd_1,\dd_2\right)=(v,w,\d_1,\d_2)\right\},
\end{equation}
where $a$ is a score in the hidden space. We have, for example
\begin{align*}
E\left\{a(X,Y)\bigm|\left(T,U,\dd_1,\dd_2\right)=(v,w,1,1)\right\}
=\frac{\int_{x\le v,\,y\le w}a(x,y)\,dF_0(x,y)}{F_0(v,w)}
\end{align*}
With this notation, we want to solve the equation
\begin{align}
\label{score_equation}
&E\left\{\th_{F_0}(T,U,\dd_1,\dd_2)\bigm|(X,Y)=(x,y)\right\}\nonumber\\
&=\int_{v\ge x,\,w\ge y}\th_{F_0}(v,w,1,1)\,g(v,w)\,dv\,dw
+\int_{v\ge x,\,w<y}\th_{F_0}(v,w,1,0)\,g(v,w)\,dv\,dw\nonumber\\
&\quad+\int_{v<x,\,w\ge y}\th_{F_0}(v,w,0,1)\,g(v,w)\,dv\,dw
+\int_{v<x,\,w<y}\th_{F_0}(v,w,0,0)\,g(v,w)\,dv\,dw\nonumber\\
&=\k_{(t,u)}(x,y).
\end{align}
Defining
\begin{equation}
\label{definition_phi}
\f_{F_0}(x,y)=\int_{v\le x,\,w\le y}a(v,w)\,dF_0(v,w),
\end{equation}
where $a$ is a score function in the hidden space,
and differentiating (\ref{score_equation}) w.r.t.\ $x$ and $y$, we now obtain the equation:
\begin{align*}
&\frac{\f_{F_0}(x,y)}{F_0(x,y)}-\frac{\f_{F_0}(x,1)-\f_{F_0}(x,y)}{F_0(x,1)-F_0(x,y)}
-\frac{\f_{F_0}(1,y)-\f_{F_0}(x,y)}{F_0(1,y)-F_0(x,y)}
-\frac{\f_{F_0}(x,1)+\f_{F_0}(1,y)-\f_{F_0}(x,y)}{1-F_0(x,1)-F_0(1,y)+F_0(x,y)}\nonumber\\
&=g(x,y)^{-1}\frac{\partial^2\k_{(t,u)}(x,y)}{\partial x\partial y}=\frac{\left\{K_h(t-x)+K_h(2-t-x)\right\}
\left\{K_h(u-y)+K_h(2-u-y)\right\}}{g(x,y)}
\end{align*}
This equation has the solution
\begin{align}
&\f_{F_0}(x,y)\nonumber\\
&=\frac{\left\{K_h(t-x)+K_h(2-t-x)\right\}
\left\{K_h(u-y)+K_h(2-u-y)\right\}}{g(x,y)}\nonumber\\
&\quad\cdot\left\{\frac{1}{F_0(x,y)}+\frac{1}{F_0(x,1)-F_0(x,y)}
+\frac{1}{F_0(1,y)-F_0(x,y)}
+\frac{1}{1-F_0(x,1)-F_0(1,y)+F_0(x,y)}\right\}^{-1}
\end{align}
Note that the solution satisfies:
$$
\f_{F_0}(1,y)=\f_{F_0}(x,1)=0,\,x,y\in[0,1].
$$

This suggests that the asymptotic behavior of the SMLE is given by:
\begin{align}
\label{asymp_representation}
\int \th_{F_0}(x,y,\d_1,\d_2)\,d\left(\P_n-P\right)(x,y,\d_1,\d_2),
\end{align}
where
\begin{align*}
&\th_{F_0}(v,w,\d_1,\d_2)\\
&=\frac{\d_1\d_2\f_{F_0}(x,y)}{F_0(x,y)}-\frac{\d_1(1-\d_2)\f_{F_0}(x,y)}{F_0(x,1)-F_0(x,y)}
-\frac{(1-\d_1)\d_2\f_{F_0}(x,y)}{F_0(1,y)-F_0(x,y)}
+\frac{(1-\d_1)(1-\d_2)\f_{F_0}(x,y)}{1-F_0(1,y)-F_0(x,1)+F_0(x,y)}\,,
\end{align*}
leading at interior points $(t,u)$ to an asymptotic variance, given by:
\begin{align*}
&\frac1n\int_{(x,y)\in[0,1]^2} \left\{\frac1{F_0(x,y)}+\frac1{F_0(x,1)-F_0(x,y)}
+\frac1{F_0(1,y)-F_0(x,y)}
+\frac1{1-F_0(1,y)-F_0(x,1)+F_0(x,y)}\right\}^{-1}\\
&\qquad\qquad\qquad\qquad\qquad\qquad\cdot \frac{\left\{K_h(t-x)+K_h(2-t-x)\right\}^2
\left\{K_h(u-y)+K_h(2-u-y)\right\}^2}{g(x,y)}\,dx\,dy\\
&\sim \left(nh^2\right)^{-1}\left\{\frac1{F_0(t,u)}+\frac1{F_0(t,1)-F_0(t,u)}
+\frac1{F_0(1,t)-F_0(t,u)}
+\frac1{1-F_0(1,u)-F_0(t,1)+F_0(t,u)}\right\}^{-1}\\
&\qquad\qquad\qquad\qquad\qquad\qquad\qquad\qquad\qquad\qquad\qquad\qquad\qquad\qquad\qquad\qquad\qquad\cdot g(t,u)^{-1}\left\{\int K(v)^2\,dv\right\}^2.
\end{align*}

Assume that $\max(t,u)\le 1-h$. Then the bias is given by:
\begin{align*}
&\int \IK_h(t-v)\IK_h(u-w)f_0(v,w)\,dv\,dw-F_0(t,u)\\
&=\int \left\{\int K_h(t-v)\int_0^v f_0(x,w)\,dx\,dv\right\}\IK_h(u-w)\,dw-F_0(t,u)\\
&=\int K_h(t-v)K_h(u-w)F_0(v,w)\,dv\,dw-F_0(t,u)\\
&=\int K(v)K(w)\left\{F_0(t-hv,u-hw)-F_0(t,u)\right\}\,dv\,dw\\
&=\tfrac12\left\{\partial_1^2F_0(t,u)+\partial_2^2F_0(t,u)\right\}h^2\left\{\int x^2 K(x)\,dx\right\}^2+o\left(h^2\right).
\end{align*}

\begin{figure}[!ht]
\begin{subfigure}[b]{0.45\textwidth}
\centering
\includegraphics[width=\textwidth]{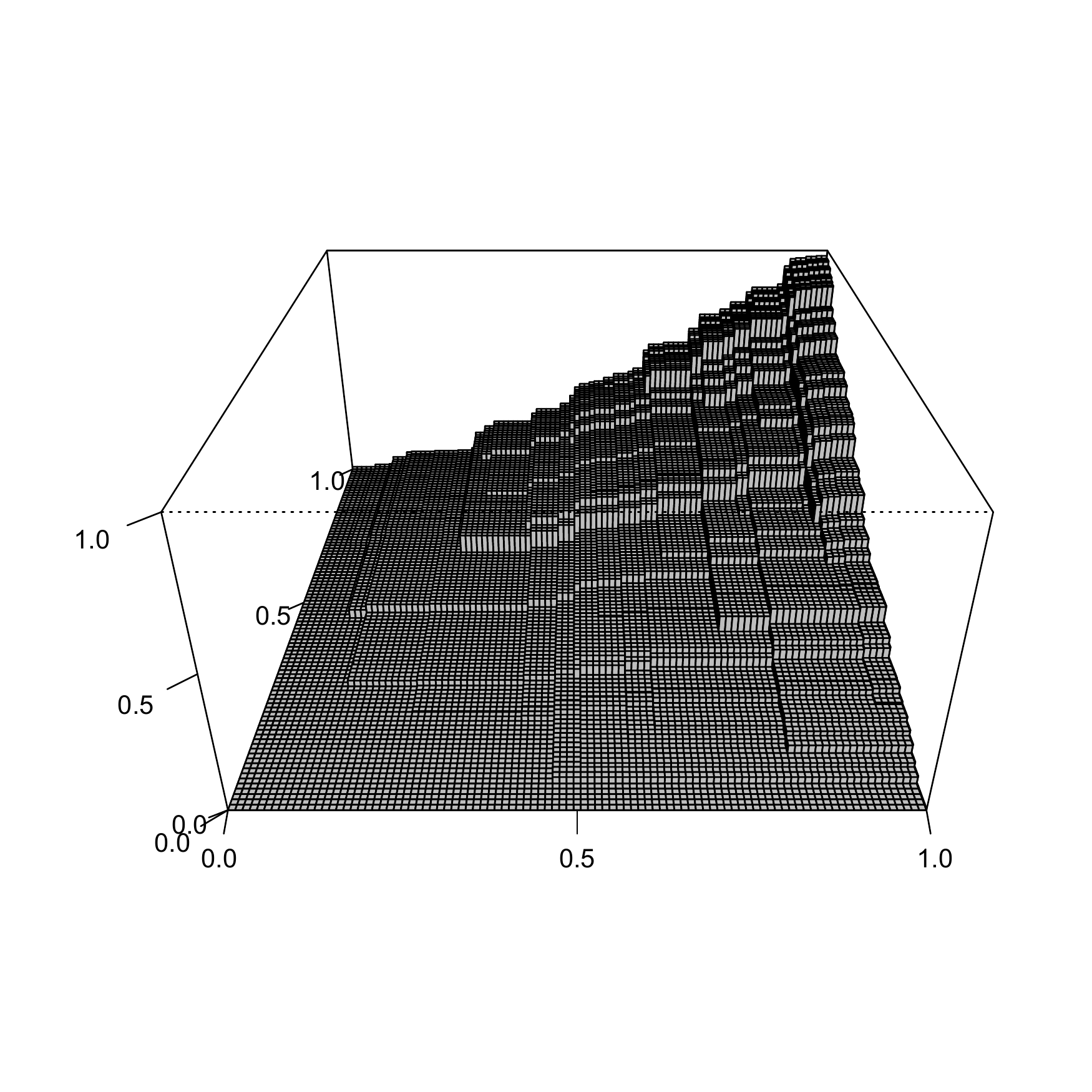}
\caption{MLE}
\label{fig:MLE_bivar2}
\end{subfigure}
\begin{subfigure}[b]{0.45\textwidth}
\centering
\includegraphics[width=\textwidth]{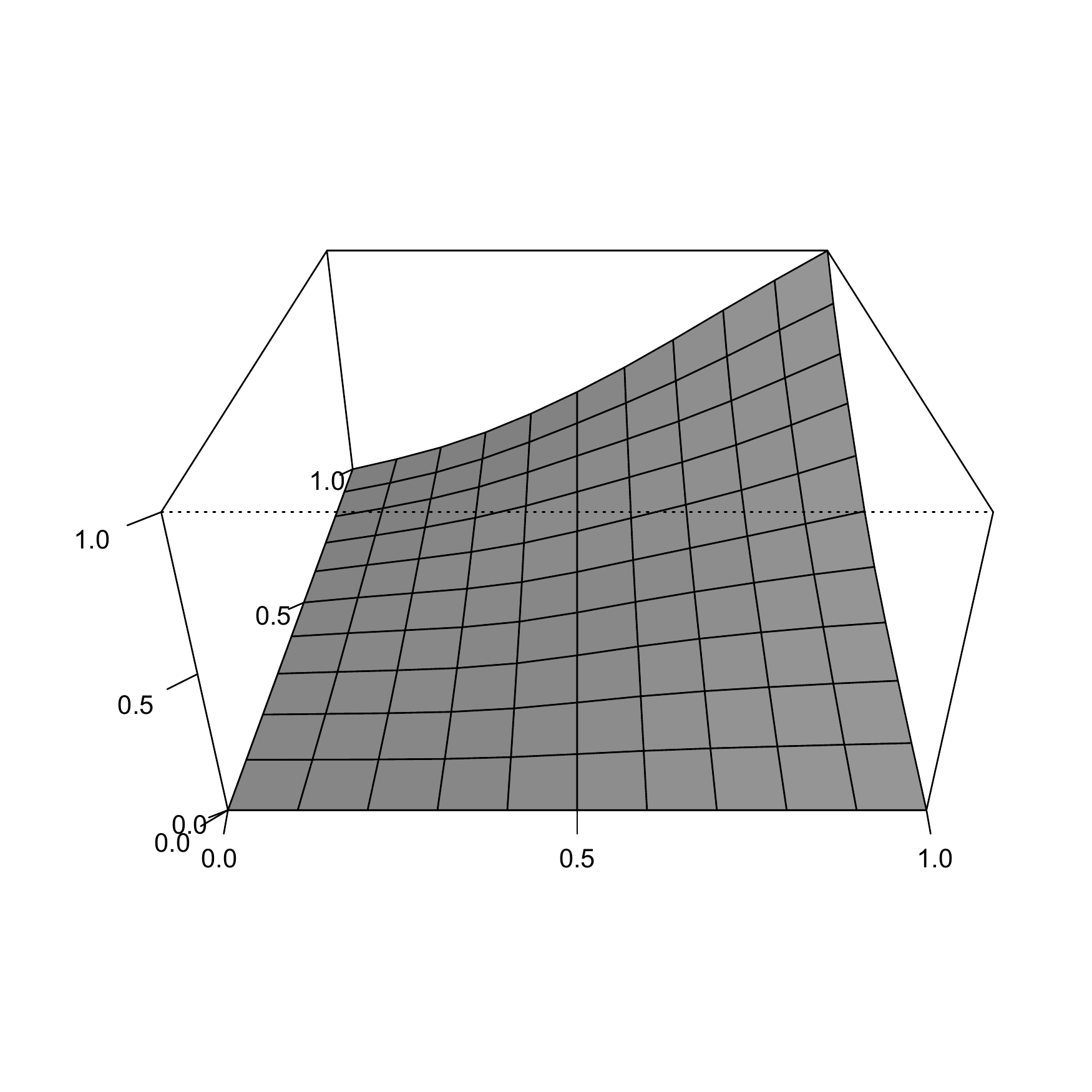}
\caption{SMLE}
\label{fig:SMLE_bivar}
\end{subfigure}
\caption{The MLE and SMLE for for a sample of size $n=1000$ of bivariate current status data, where the hidden variables have a distribution with density $f_0(x,y)=x+y$, and the observation distribution is uniform on $[0,1]^2$.}
\label{fig:MLE+SMLE_bivarCS}
\end{figure}

The SMLE is compared with the MLE in Figure \ref{fig:MLE+SMLE_bivarCS}.
Using the assumption that $\hat F_n^{(SML)}(t,u)-F_0(t,u)$ has the asymptotic representation (\ref{asymp_representation}), we get the following result.

\begin{theorem}[Conjecture]
\label{th:SMLE}
Under the conditions of Theorem \ref{th:plug-in} we have for each point $(t,u)\in(0,1)^2$ satisfying these conditions, if $c=\lim_{n\to\infty}n^{1/3}h_n^2$,
$$
n^{1/3}\left\{\hat F_{n,h_n}^{(SML)}(t,u)-F_0(t,u)\right\}\stackrel{{\cal D}}\longrightarrow
N\left(\beta,\s^2\right),
$$
where $N(\b,\s^2)$ is a normal distribution with first moment
$$
\b=\tfrac12c\left\{\partial_1^2F_0(t,u)+\partial_2^2F_0(t,u)\right\}\left\{\int x^2 K(x)\,dx\right\}^2,
$$
and variance
\begin{align*}
\s^2&=c^{-1}\left\{\frac1{F_0(t,u)}+\frac1{F_0(t,1)-F_0(t,u)}
+\frac1{F_0(1,t)-F_0(t,u)}
+\frac1{1-F_0(1,u)-F_0(t,1)+F_0(t,u)}\right\}^{-1}\\
&\qquad\qquad\qquad\qquad\qquad\qquad\qquad\qquad\qquad\qquad\qquad\qquad\qquad\qquad\qquad\cdot g(t,u)^{-1}\left\{\int K(v)^2\,dv\right\}^2.
\end{align*}
\end{theorem}

\begin{remark}
{\rm
Note that choosing $h_n\asymp n^{-1/6}$ is the asymptotically optimal choice (modulo constants) since the variance is of order $1/(nh_n^2)$ and the bias of order $h_n^2$, unless the bias is of order $o(h_n^2)$ (as happens when $F_0$ is the uniform distribution function on $[0,1]^2$). Also note that the bias term, caused by the interaction of the observation distribution $G$ and the distribution function $F_0$, which entered into the bias term of Theorem \ref{th:plug-in}, plays no role here.
}
\end{remark}

\section{A simulation study}
\label{section:simulations}
In order to compare the behavior of the three estimators, we took 1000 samples from the distribution with distribution function
$$
F_0(x,y)=\tfrac12xy(x+y),\,x,y\in[0,1]^2,
$$
and generated bivariate current status data from this with respect to the uniform distribution on $[0,1]^2$. The samples size taken were $n=100, 500, 1000$ and $5000$, respectively. So we have observations $(T_i,U_i)$ from the uniform distribution in $[0,1]^2$, and for each such pair we get from the corresponding pair $(X_i,Y_i)$, drawn from $F_0$, independently w.r.t.\ $(T_i,U_i)$, the indicators
$$
\dd_{i1}=1_{\{X_i\le T_i\}}\mbox{ and }\dd_{i2}=1_{\{Y_i\le U_i\}}.
$$

\begin{table}[!ht]
\caption{Estimated values of the standard deviations times $n^{1/3}$ for three estimators of $F_0(t,u)$ at a number of values of $(t,u)$, where $F_0(t,u)=\tfrac12tu(t+u)$. The values corresponding to $\infty$ are the asymptotic values, deduced from Theorems \ref{th:plug-in} and \ref{th:SMLE}. The asymptotic values for the MLE are unknown.}
\label{table:simstudy}
\begin{subfigure}[b]{0.3\textwidth}
\centering
\begin{tabular}{|c|c|c|c|c|}
\hline
 & \multicolumn{4}{c|}{$u=0.6$}\\
\hline
$t$ &  $\mbox{$n$}$ & MLE &  SMLE & Plug-in\\
\hline
$0.2$ & $100$ & $0.251$ & $0.133$ & $0.185$ \\
 & $500$ & $0.257 $ & $0.132$ & $0.152$  \\
 & $1000$ & $0.247$ & $0.128$ & $0.142$\\
 & $5000$ & $0.246$ & $0.123$ & $0.126$ \\
  & $\infty$ & $-$ & $0.130$ & $0.107$ \\
\hline
$0.4$ & $100$ & $0.379$ & $0.212$ & $0.198$ \\
 & $500$ & $0.367$ & $0.194$ & $0.176$ \\
 &  $1000$ & $0.360$ & $0.190$ & $0.179$ \\
 & $5000$ & $0.357$ & $0.178$ & $0.156$ \\
   & $\infty$ & $-$ & $0.182$ & $0.162$ \\
\hline
$0.6$ & $100$ & $ 0.475$ &  $0.263$ & $0.226$ \\
 & $500$ & $0.412$ & $0.223$ & $0.221$ \\
 &  $1000$ & $0.450$ & $0.238$ & $0.215$ \\
 &  $5000$ & $0.441$ & $0.218$ & $0.205$ \\
   & $\infty$ & $-$ & $0.203$ & $0.206$ \\
\hline
$0.8$ & $100$ & $0.550$ &  $0.276$ & $0.290$ \\
 & $500$ & $0.508$ & $0.253$ & $0.265$ \\
 &  $1000$ & $0.503$ & $0.255$ & $0.266$ \\
 & $5000$ & $0.514$ & $0.236$ & $0.240$\\
   & $\infty$ & $-$ & $0.183$ & $0.236$ \\
\hline
\end{tabular}
\caption{$n^{1/3}$ times standard deviation}
\label{tab:sd}
\end{subfigure}
\hspace{3cm}
\begin{subfigure}[b]{0.3\textwidth}
\centering
\begin{tabular}{|c|r|r|c|c|}
\hline
 & \multicolumn{4}{c|}{$u=0.6$}\\
\hline
$t$ &  $\mbox{$n$}$ & MLE &  SMLE & Plug-in\\
\hline
$0.2$ & $100$ & $0.003$ & $0.043$ & $0.205$ \\
 & $500$ & $-0.001 $ & $0.037$ & $0.263$  \\
 & $1000$ & $-0.037$ & $0.029$ & $0.264$\\
 & $5000$ & $\!\!\!-0.0002$ & $0.030$ & $0.230$ \\
  & $\infty$ & $-$ & $0.044$ & $0.133$ \\
\hline
$0.4$ & $100$ & $\!\!\!\!\!-0.004$ & $0.061$ & $0.155$ \\
 & $500$ & $-0.006$ & $0.056$ & $0.169$ \\
 &  $1000$ & $-0.022$ & $0.048$ & $0.164$ \\
 & $5000$ & $\!\!\!-0.006$ & $0.049$ & $0.167$ \\
   & $\infty$ & $-$ & $0.056$ & $0.166$ \\
\hline
$0.6$ & $100$ & $0.115$ &  $0.099$ & $0.015$ \\
 & $500$ & $0.003$ & $0.084$ & $0.202$ \\
 &  $1000$ & $-0.036$ & $0.080$ & $0.190$ \\
 &  $5000$ & $0.007$ & $0.072$ & $0.188$ \\
   & $\infty$ & $-$ & $0.067$ & $0.200$ \\
\hline
$0.8$ & $100$ & $0.118$ &  $0.124$ & $-0.152$ \\
 & $500$ & $-0.018$ & $0.112$ & $-0.078$ \\
 &  $1000$ & $0.006$ & $0.117$ & $-0.090$ \\
 & $5000$ & $-0.0002$ & $0.116$ & $-0.006$\\
   & $\infty$ & $-$ & $0.078$ & \,\,\,\,\,$0.233$ \\
\hline
\end{tabular}
\caption{$n^{1/3}$ times bias}
\label{tab:bias}
\end{subfigure}
\end{table}

Since simulations for sample size $5000$ with the `real' MLE, based on the maximal intersection rectangles, would have taken prohibitively long, we instead computed the MLE on a sieve, constructed in the following way. For each sample of size $n$, we distributed $m_n=\left[n^{2/3}\right]$ points on the unit square, by letting their $x$- and $y$-coordinates be multiples of $n^{-1/3}$ and  permuting these coordinates randomly, according to the uniform distribution on permutations. Here $[x]$ denotes the largest integer $\le x$.

\begin{figure}[!ht]
\begin{subfigure}[b]{0.45\textwidth}
\centering
\includegraphics[width=\textwidth]{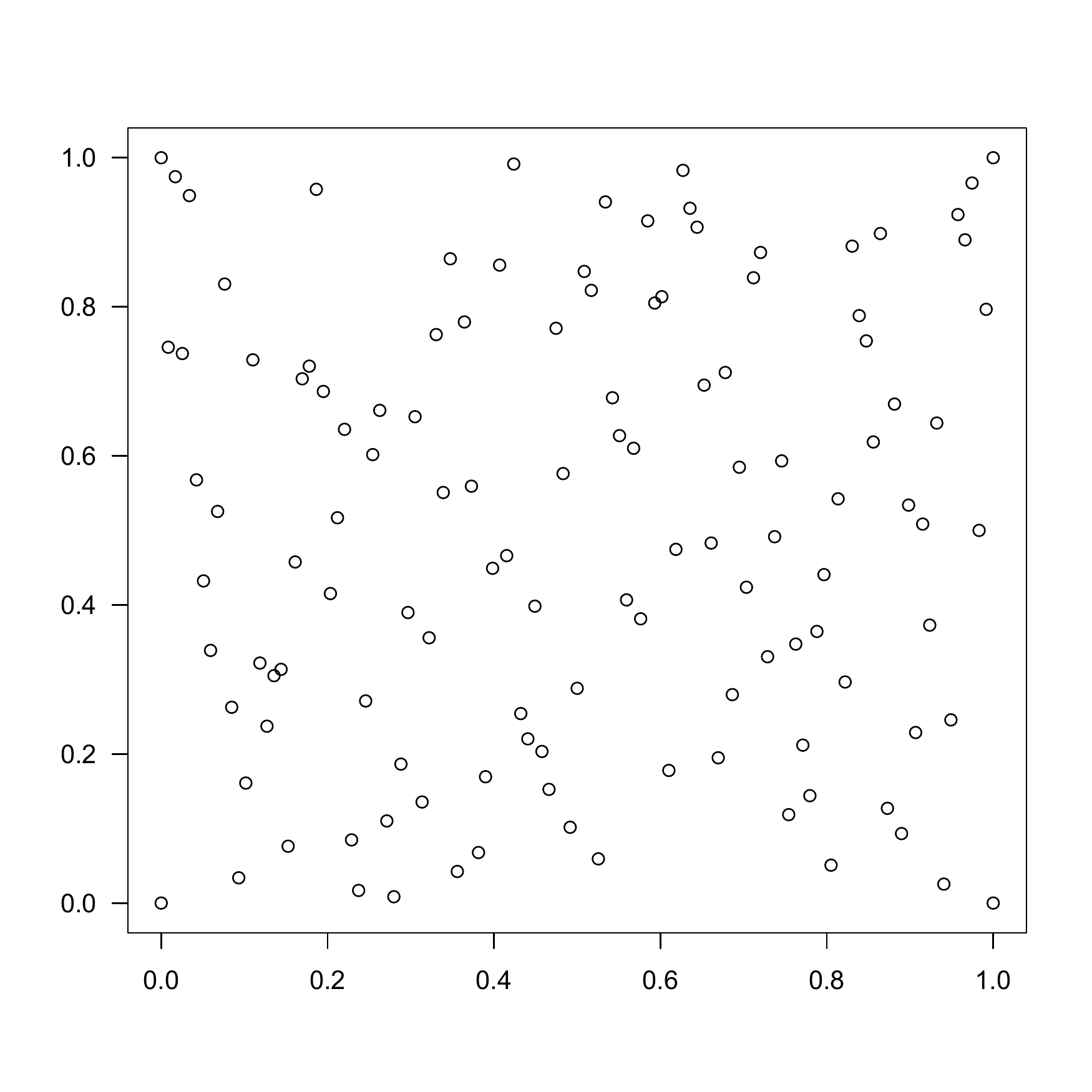}
\caption{Initial set of points of possible mass}
\label{fig:initial_mass}
\end{subfigure}
\begin{subfigure}[b]{0.45\textwidth}
\centering
\includegraphics[width=\textwidth]{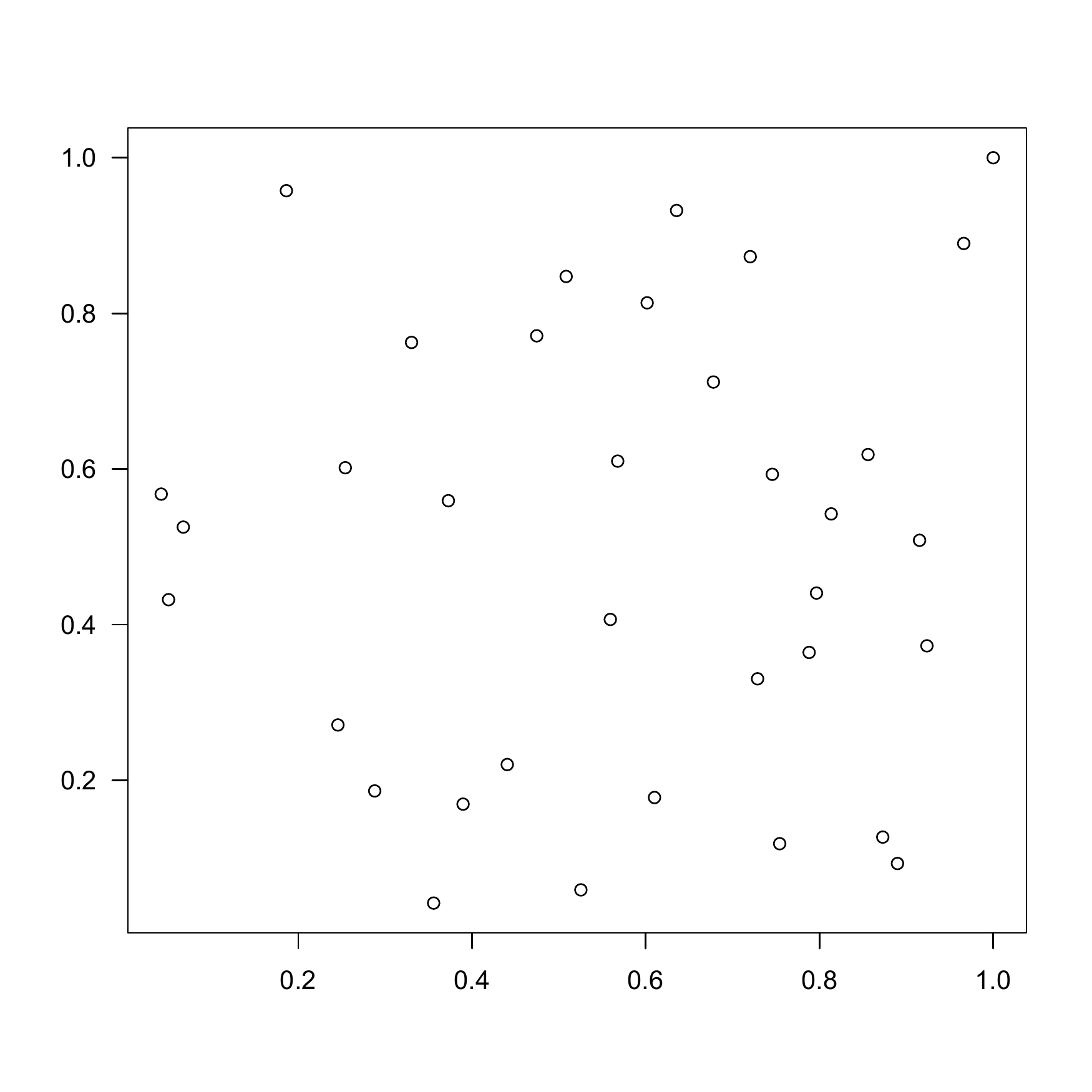}
\caption{Mass points of MLE}
\label{fig:MLE_mass}
\end{subfigure}
\caption{A set of points on which the MLE is allowed to put its mass and the actual set of mass points of the MLE for a sample of size $n=1000$.}
\label{fig:massMLE}
\end{figure}

So we start with order $n^{2/3}$ points on which the sieved MLE can place its mass, to which we add the vertices of the unit square to ensure a finite log likelihood, and compute the MLE which only is allowed to have mass at these points. This set of points is shown in Figure \ref{fig:massMLE}, corresponding to sample size $n=1000$, and the reduced set of points on which the MLE actually puts positive mass is also shown in this picture.

\begin{figure}[!ht]
\centering
\includegraphics[width=0.5\textwidth]{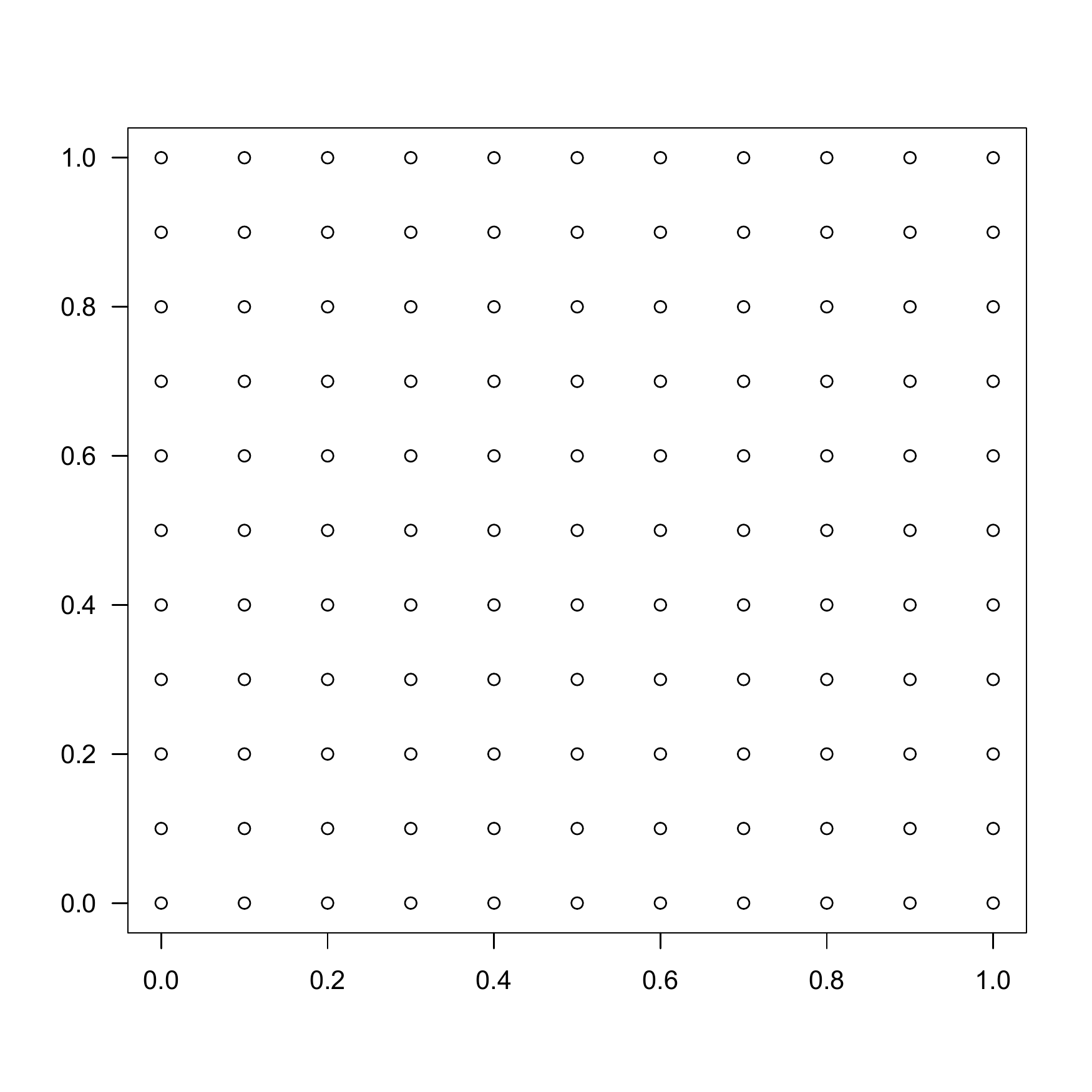}
\caption{The set of points on which the plug-in estimator puts its mass for $n=1000$.}
\label{fig:massPlugin}
\end{figure}

The rather different set of points on which the plug-in estimate puts its mass is shown in Figure \ref{fig:massPlugin}. In fact, inspection of the set of points on which the MLE, based on the maximal intersection rectangles puts its mass, shows that that such sets are rather similar to Figure \ref{fig:MLE_mass} and not similar to Figure \ref{fig:massPlugin}. By the rather irregular structure of sets like Figure \ref{fig:MLE_mass} the bias of the MLE is reduced w.r.t.\ the plug-in estimator which is constant on squares with sides of order $n^{-1/3}$. We note, however, that the number of points in Figure \ref{fig:initial_mass} is the same as in Figure \ref{fig:massPlugin} (the number is $121$).

\begin{table}[!h]
\caption{Estimated values of the standard deviations times $n^{1/3}$ for three estimators of $F_0(t,u)$ at a number of values of $(t,u)$, where $F_0(t,u)=tu$. The values corresponding to $\infty$ are the asymptotic values, deduced from Theorems \ref{th:plug-in} and \ref{th:SMLE}. The asymptotic values for the MLE are unknown.}
\label{table:simstudy2}
\begin{subfigure}[b]{0.3\textwidth}
\centering
\begin{tabular}{|c|c|c|c|c|}
\hline
 & \multicolumn{4}{c|}{$u=0.6$}\\
\hline
$t$ &  $\mbox{$n$}$ & MLE &  SMLE & Plug-in\\
\hline
$0.2$ & $100$ & $0.371$ & $0.199$ & $0.388$ \\
 & $500$ & $0.367 $ & $0.173$ & $0.286$  \\
 & $1000$ & $0.363$ & $0.145$ & $0.262$\\
 & $5000$ & $0.360$ & $0.121$ & $0.195$ \\
  & $\infty$ & $-$ & $0$ & $0$ \\
\hline
$0.4$ & $100$ & $0.458$ & $0.266$ & $0.511$ \\
 & $500$ & $0.438$ & $0.220$ & $0.391$ \\
 &  $1000$ & $0.447$ & $0.197$ & $0.332$ \\
 & $5000$ & $0.438$ & $0.159$ & $0.257$ \\
   & $\infty$ & $-$ & $0$ & $0$ \\
\hline
$0.6$ & $100$ & $ 0.480$ &  $0.294$ & $0.576$ \\
 & $500$ & $0.486$ & $0.239$ & $0.433$ \\
 &  $1000$ & $0.491$ & $0.221$ & $0.362$ \\
 &  $5000$ & $0.470$ & $0.178$ & $0.301$ \\
   & $\infty$ & $-$ & $0$ & $0$ \\
\hline
$0.8$ & $100$ & $0.517$ &  $0.303$ & $0.602$ \\
 & $500$ & $0.497$ & $0.242$ & $0.438$ \\
 &  $1000$ & $0.490$ & $0.226$ & $0.386$ \\
 & $5000$ & $0.479$ & $0.177$ & $0.303$\\
   & $\infty$ & $-$ & $0$ & $0$ \\
\hline
\end{tabular}
\caption{$n^{1/3}$ times standard deviation}
\label{tab:sd2}
\end{subfigure}
\hspace{3cm}
\begin{subfigure}[b]{0.3\textwidth}
\centering
\begin{tabular}{|c|r|r|r|r|}
\hline
 & \multicolumn{4}{c|}{$u=0.6$}\\
\hline
$t$ &  $\mbox{$n$}$ & MLE &  SMLE & Plug-in\\
\hline
$0.2$ & $100$ & $-0.032$ & $-0.009$ & $-0.007$ \\
 & $500$ & $0.001 $ & $0.037$ & $0.005$  \\
 & $1000$ & $-0.027$ & $-0.006$ & $-0.008$\\
 & $5000$ & $0.006$ & $0.024$ & $-0.009$ \\
  & $\infty$ & $-$ & $0$ & $0$ \\
\hline
$0.4$ & $100$ & $-0.016$ & $-0.009$ & $-0.028$ \\
 & $500$ & $0.018$ & $0.007$ & $0.004$ \\
 &  $1000$ & $-0.022$ & $-0.001$ & $0.001$ \\
 & $5000$ & $-0.008$ & $-0.003$ & $0.004$ \\
   & $\infty$ & $-$ & $0$ & $0$ \\
\hline
$0.6$ & $100$ & $0.106$ &  $0.023$ & $-0.025$ \\
 & $500$ & $0.011$ & $0.022$ & $0.014$ \\
 &  $1000$ & $-0.004$ & $0.005$ & $-0.008$ \\
 &  $5000$ & $0.027$ & $0.010$ & $0.004$ \\
   & $\infty$ & $-$ & $0$ & $0$ \\
\hline
$0.8$ & $100$ & $0.092$ &  $0.033$ & $-0.011$ \\
 & $500$ & $-0.015$ & $0.037$ & $-0.001$ \\
 &  $1000$ & $0.012$ & $0.012$ & $0.008$ \\
 & $5000$ & $0.014$ & $-0.011$ & $0.006$\\
   & $\infty$ & $-$ & $0$ &$0$ \\
\hline
\end{tabular}
\caption{$n^{1/3}$ times bias}
\label{tab:bias2}
\end{subfigure}
\end{table}

We took the bandwidth $h_n$ for the SMLE equal to $n^{-1/6}$ and we also used this as binwidth for the plug-in estimator.
Table \ref{tab:sd} shows that there is no indication that $n^{1/3}$ times the standard deviation of the MLE is increasing with sample size, and Table \ref{tab:bias} suggests that the bias times $n^{1/3}$ is vanishing (as is also true in the one-dimensional case!), as $n\to\infty$, so the hypothesis that the rate of convergence of the MLE is $n^{1/3}$ is not contradicted.

However, if the rate of the MLE is actually of order $n^{1/3}(\log n)^{-1/3}$, for example, we can probably not detect this in the present way. The theory for the MSLE and plug-in estimators seems to be confirmed by the simulations, in particular at the points $(0.4,0.6)$ and $(0.6,0,6)$, where the boundary effects are still not active. Note that if, for example, $n=500$, the bandwidths for the SMLE are equal to $500^{-1/6}\approx0.3549537$, so the boundary kernel starts getting active for the points $(0.2,0.6)$ and $(0.8,0.6)$. This is still true for sample size $n=5000$, where the bandwidth is $5000^{-1/6}\approx0.2418271$. It is clear from Table \ref{tab:sd} that the MLE has a bigger variance than the other two estimators, but on the other hand the bias of the MLE is usually smaller than that of the other estimators.

If the second order partial derivatives of the distribution function $F_0$ vanish at $(t,u)$, as happens, for example, with the uniform distribution, it is possible to achieve higher rates of convergence with the SMLE and the plug-in estimator. We demonstrate this for the uniform distribution function $F_0$, where we keep the bandwidth constant and equal to $0.4$ for the SMLE and equal to 0.2 for the plug-in estimator (to avoid the boundary correction). In this case the bias times $n^{1/3}$ tends to zero for the SMLE and the plug-in estimator, see Table \ref{tab:bias2}. If we keep the bandwidth constant, the variance of the SMLE and plug-in estimator should be of order $n^{-1}$ and these estimators should therefore actually attain a parametric rate of convergence in this case. We expect the MLE to have again rate $n^{1/3}$ in this case however, which is also (to a certain extent) suggested by Table \ref{table:simstudy2}.

\section{More general bivariate interval censoring}
\label{section:interval_cens}

The purpose of this section is to show that the MLE and SMLE, discussed in the preceding sections in the context of the bivariate current status model, can also be used for more general interval censored data. As mentioned earlier, the current status model is the simplest case of the interval censoring model. For the bivariate interval censoring, case 2, model, the data are of the form
$$
\left(T_{i1},U_{i1},T_{i2},U_{i2},\dd_{i1}^{(1)},\dd_{i2}^{(1)},\dd_{i1}^{(2)},\dd_{i2}^{(2)}\right),\,i=1,\dots,n,
$$
where
$$
\dd_{i1}^{(1)}=1_{\{X_{i1}\le T_{i1}\}},\quad\dd_{i1}^{(2)}=1_{\{T_{i1}<X_{i1}\le U_{i1}\}},
\quad\dd_{i2}^{(1)}=1_{\{X_{i2}\le T_{i2}\}},\quad\dd_{i2}^{(2)}=1_{\{T_{i2}<X_{i2}\le U_{i2}\}},
$$
Defining
$$
\dd_{i1}^{(3)}=1-\dd_{i1}^{(1)}-\dd_{i1}^{(2)},\qquad\dd_{i2}^{(3)}=1-\dd_{i2}^{(1)}-\dd_{i2}^{(2)},
$$
and defining the corresponding generic values $\d_i^{(j)}$ similarly, we can define  the measure $dV_n^{(ij)}$ by
$$
dV_n^{(ij)}=\d_1^{(i)}\d_1^{(j)}\d_2^{(i)}\d_2^{(j)}\,d\P_n\left(t,u,v,w,\d_1^{(i)},\d_2^{(i)},\d_1^{(j)},\d_2^{(j)}\right),\,1\le i,j\le3,
$$
and an MLE of $F$ is then obtained by maximizing
\begin{align}
\label{bivarICM_likelihood}
\ell(F)&=\int \log F(t,v)\,dV_n^{(11)}
+\int_{t\le u,\,v\le w}\log\left\{F(u,w)-F(t,w)-F(u,v)+F(t,v)\right\}\,dV_n^{(22)}\nonumber\\
&\qquad+\int_{u\ge t}\log\left\{F(u,v)-F(t,v)\right\}
\,dV_n^{(21)}+\int_{w\ge v} \log\left\{F(t,w)-F(t,v)\right\}
\,dV_n^{(12)}\nonumber\\
&\qquad+\int_{u\ge t}\log\left\{F_1(u)-F_1(t)-F(u,w)+F(t,w)\right\}
\,dV_n^{(23)}\nonumber\\
&\qquad+\int_{w\ge v}\log\left\{F_2(w)-F_2(v)-F(u,w)+F(u,v)\right\}
\,dV_n^{(32)}\nonumber\\
&\qquad+\int_{w\ge v} \log\left\{F_1(t)-F(t,w)\right\}
\,dV_n^{(13)}+\int_{u\ge t}\log\left\{F_2(v)-F(u,v)\right\}
\,dV_n^{(31)}\nonumber\\
&\qquad+\int \log\left\{1-F_1(u)-F_2(w)+F(u,w)\right\}\,dV_n^{(33)}
\end{align}
over bivariate distribution functions $F$.
The Fenchel duality conditions become:
\begin{align}
\label{fenchel-ineq-ML2}
&\int_{t\ge x,\,v\ge y} \frac{dV_n^{(11)}}{F(t,v)}
+\int_{t<x\le u,\,v<y\le w}
\frac{dV_n^{(22)}}{F(u,w)-F(t,w)-F(u,v)+F(t,v)}\nonumber\\
&\qquad+\int_{t<x\le u,\,v\ge y} \frac{dV_n^{(21)}}{F(u,v)-F(t,v)}
+\int_{t\ge x,\,v<y\le w} \frac{dV_n^{(12)}}{F(t,w)-F(t,v)}\nonumber\\
&\qquad+\int_{t<x\le u,\,w\ge y}
\frac{dV_n^{(23)}}{F_1(u)-F_1(t)-F(u,w)+F(t,w)}\nonumber\\
&\qquad+\int_{u\ge x,\,v<y\le w}\frac{dV_n^{(32)}}{F_2(w)-F_2(v)-F(u,w)+F(u,v)}\nonumber\\
&\qquad+\int_{t<x\le u,\,v\ge y} \frac{dV_n^{(13)}}{F_1(t)-F(t,w)}+\int_{t\ge x,\,v<y\le w} \frac{dV_n^{(31)}}{F_2(v)-F(u,v)}\nonumber\\
&\qquad+\int_{u<x,\,w<y} \frac{dV_n^{(33)}}{1-F_1(u)-F_2(w)+F(u,w)}\le1,
\end{align}
with equality if $(x,y)$ is a point of mass of $dF$. 

For computational purposes (but probably not for the development of distribution theory!) it is more convenient not to distinguish between the measures $V_n^{(ij)}$ and just to introduce rectangles to which the unobservable observations are known to belong, where we represent the (one-sided) unbounded rectangles by finite rectangles with upper or lower bounds outside the range of the observed data. In this set-up we simply have to maximize
\begin{equation}
\label{simple_loglike_bivarIC}
\sum f_i\log H_i'p,
\end{equation}
where $p=(p_1,\dots,p_m)'$ is a vector of probability masses at possible points of mass $(x_j,y_j)$ and $H_i$ is a vector of length $m$, consisting of ones and zeros, where the component $H_{ij}$ is equal to $1$ if the point $(x_j,y_j)$ is contained in the rectangle
$$
[L_{i1},R_{i1}]\times[L_{i2},R_{i2}],
$$
and is zero, otherwise, and where the $f_i$ denote the multiplicities at the $i$th observation point. The masses $p_j$ should be nonnegative and sum to $1$. This optimization can easily be handled by using iterative quadratic minimization and the support reduction algorithm, documented in \cite{piet_geurt_jon:08}. In fact, the treatment is completely analogous to the treatment of the Aspect experiment for quantum statistics, discussed there.

The data of \cite{Betensky:99} are given in Table \ref{tab:Betensky}, where the rectangles to which the hidden observations are known to belong are denoted by $[L_{i1},R_{i1}]\times[L_{i2},R_{i2}]$, $i=1,\dots,n$. The frequencies of the hidden observations belonging to these rectangles are given in the $5$th and $10$th column. There are $87$ observation rectangles and the total sample size, taking the frequencies into account, is $204$. The table is also given in \cite{sun:06}, Table 7.1, p.\ 165, but there the rectangles are slightly changed from the data in \cite{Betensky:99} by lowering the left bounds by $1$ if they are larger than zero. Since we do not see a pressing reason for doing that, we just give the data here as they were given by \cite{Betensky:99}. If the upper bound $R_{ij}$ is unknown, we put $R_{ij}=\infty$ and if the lower bound $L_{ij}$ is unknown, we put $L_{ij}=-\infty$.

\begin{table}[!ht]
\begin{centering}
\caption{The Betensky-Finkelstein data}
\label{tab:Betensky}
\begin{tabular}{|c|r|r|r|c|r|r|r|r|r|c|}
\hline
$L_{i1}$ &$R_{i1}$ &$L_{i2}$ &$R_{i2}$ &frequency  &\qquad\qquad\qquad   &$L_{i1}$ &$R_{i1}$ &$L_{i2}$ &$R_{i2}$ &frequency\\
\hline
  0  &  3  &  0  &$-$  &   3 &  & 6  &$-$  &  6  &$-$  &   3\\
  0  &  3  &  3  &$-$  &   1 &  & 6  &$-$  &  9  &$-$  &   1\\
  0  &  3  &  6  &$-$  &   3 &  & 9  &$-$  &  0  &$-$  &   2\\
  0  &  6  &  6  &$-$  &   1 &  & 9  &$-$  &  9  &$-$  &   3\\
  0  &  3  &  9  &$-$  &   1 &  & 9  &$-$  & 12  &$-$  &   1\\
  0  &  3  & 12  &$-$  &   5 & & 12  &$-$  &  0  &$-$  &   5\\
  0  &  3  & 15  &$-$  &   5 & & 12  &$-$  &  6  &$-$  &   1\\
  0  &  6  & 15  &$-$  &   1 & & 12  &$-$  &  9  &$-$  &   4\\
  3  &  3  &  3  &$-$  &   1 & & 12  &$-$  & 12  &$-$  &  10\\
  3  &  3  &  6  &$-$  &   1 & & 15  &$-$  &  0  &$-$  &   3\\
  3  &  3  &  9  &$-$  &   3 & & 15  &$-$  &  3  &$-$  &   1\\
  3  &  6  &  9  &$-$  &   2 & & 15  &$-$  &  6  &$-$  &   1\\
  3  &  6  & 12  &$-$  &   3 & & 15  &$-$  &  9  &$-$  &   2\\
  3  &  3  & 15  &$-$  &   2 & & 15  &$-$  & 12  &$-$  &   8\\
  3  &  6  & 15  &$-$  &   2 & & 15  &$-$  & 15  &$-$  &   9\\
  3  &  6  & 18  &$-$  &   1 & & 18  &$-$  &  0  &$-$  &   1\\
  3  &  3  & 21  &$-$  &   1 & & 18  &$-$  &  6  &$-$  &   1\\
  6  &  6  &  0  &$-$  &   2 & & 18  &$-$  &  9  &$-$  &   1\\
  6  &  9  &  0  &$-$  &   1 & & 18  &$-$  & 12  &$-$  &   1\\
  6  &  9  &  9  &$-$  &   1 & & 18  &$-$  & 15  &$-$  &   3\\
  6  &  6  & 12  &$-$  &   1 & & 18  &$-$  & 18  &$-$  &   6\\
  6  &  9  & 12  &$-$  &   2 & & 21  &$-$  & 15  &$-$  &   1\\
  6  &  6  & 15  &$-$  &   1 & & $-$  &  0  &  0  &$-$  &   9\\
  6  &  9  & 15  &$-$  &   1 & & $-$  &  0  &  3  &$-$  &   3\\
  6  &  6  & 18  &$-$  &   1 & & $-$  &  0  &  6  &$-$  &  10\\
  6  &  9  & 18  &$-$  &   2 & & $-$  &  0  &  9  &$-$  &   6\\
  9  &  9  &  0  &$-$  &   1 & & $-$  &  0  & 12  &$-$  &   8\\
  9  & 12  &  0  &$-$  &   2 & & $-$  &  0  & 15  &$-$  &   5\\
  9  &  9  &  9  &$-$  &   2 & & $-$  &  0  & 18  &$-$  &   4\\
  9  & 12  &  9  &$-$  &   1 & & $-$  &  0  & 21  &$-$  &   1\\
  9  & 12  & 12  &$-$  &   3 & & 0  &$-$  &  0  &  3  &   1\\
  9  &  9  & 15  &$-$  &   1 & & 6  &$-$  &  0  &  6  &   1\\
  9  & 12  & 24  &$-$  &   1 & & 6  &$-$  &  6  &  6  &   1\\
  9  &  9  & 27  &$-$  &   1 & & 12  &$-$  &  0  &  3  &   1\\
 12  & 12  &  0  &$-$  &   1 & & 12  &$-$  &  0  &  6  &   1\\
 12  & 15  &  0  &$-$  &   1 & & 15  &$-$  &  0  &  3  &   1\\
 12  & 15  &  6  &$-$  &   1 & & 21  &$-$  & 15  & 15  &   1\\
 12  & 15  & 15  &$-$  &   1 & & 3  &$-$  &$-$  &  0  &   1\\
 12  & 15  & 21  &$-$  &   1 & & 9  &$-$  &$-$  &  0  &   1\\
  0  &$-$  &  0  &$-$  &   	6 & & 12  &$-$  &$-$  &  0  &   1\\
  3  &$-$  &  0  &$-$  &   	2 & &  0  &  3  &  0  &  6  &   1\\
  6  &$-$  &  0  &$-$  &   	1 & &  3  &  6  &  6  & 12  &   1\\
  6  &$-$  &  3  &$-$  &   	2 & &  9  &  9  &  9  &  9  &   1\\
  $-$ & 0 &$-$ &0 		&1  & &  &    &   &\\
  \hline
\end{tabular}
\end{centering}
\end{table}

The maximal intersection  rectangles where the MLE will put its mass are given in Table \ref{tab:can_rectangles}. They can be computed, for example, by applying the reduction algorithm, used in the R package MLEcens.

\begin{table}[!ht]
\caption{Maximal intersection rectangles and masses of MLE}
\label{can_rectangles}
\begin{subfigure}[b]{0.3\textwidth}
\centering
\begin{tabular}{|r|r|r|r|}
\hline
$L_{j1}$ &$R_{j1}$ &$L_{j2}$ &$R_{j2}$\\
\hline
  0  &  0  &  0  &0  \\
  0  &  0  &  21  &$-$\\
  3  &  3  &  21  &$-$\\
  6  &  6  &  6  &6  \\
  6  &  6  &  18  &$-$  \\
  9  &  9  & 9  &9  \\
  9  &  9  & 27  &$-$\\
  12  &  12  & 0  &0\\
  12  &  12  &  24  &$-$\\
  15  &  15 &  0  &0\\
  15  &  15  &  21  &$-$\\
  21  &  $-$  &  15  &15\\
  21  &  $-$  & 18  &$-$\\
   \hline
\end{tabular}
\caption{Canonical rectangles}
\label{tab:can_rectangles}
\end{subfigure}
\hspace{3cm}
\begin{subfigure}[b]{0.3\textwidth}
\centering
\begin{tabular}{|r|r|c|}
\hline
$L_{j1}$ &$L_{j2}$ & mass MLE\\
\hline
  0  &  0  &  0.013676984    \\
  0  &  21  &  0.307533525  \\
  3  &  21  &  0.087051863  \\
  6  &  6  &  0.014940282    \\
  6  &  18  &  0.062521573    \\
  9  &  9  & 0.010009349    \\
  9  &  27  & 0.071073995  \\
  12  &  0  & 0.004836043  \\
  12  &  24  &  0.053334241  \\
  15  &  0 &  0.042456241  \\
  15  &  21  &  0.021573343  \\
  21  &  15  &  0.044427509 \\
  21  &  18  & 0.266565054  \\
   \hline
\end{tabular}
\caption{Masses of MLE}
\label{tab:MLE_masses}
\end{subfigure}
\end{table}

To facilitate the comparison with the existing literature, we will only discuss the MLE, based on the preliminary reduction to rectangles which can have mass, and not follow the procedure we used for computing the MLE on a sieve in the simulation from the density $f(x,y)=x+y$ on $[0,1]^2$. We will use the convention of putting the mass of the MLE in the right upper corner of these rectangles, and compute the MLE by the support reduction algorithm of \cite{piet_geurt_jon:08}. The result is shown in Table \ref{can_rectangles}, where the masses of the MLE are given. It is seen that this is in close correspondence with Table 7.2 on p.\ 166 of \cite{sun:06}, apart from the slightly different definition of the rectangles. The R package MLEcens also gives this result (in all $9$ decimals).

\begin{figure}[!ht]
\begin{subfigure}[b]{0.45\textwidth}
\centering
\includegraphics[width=\textwidth]{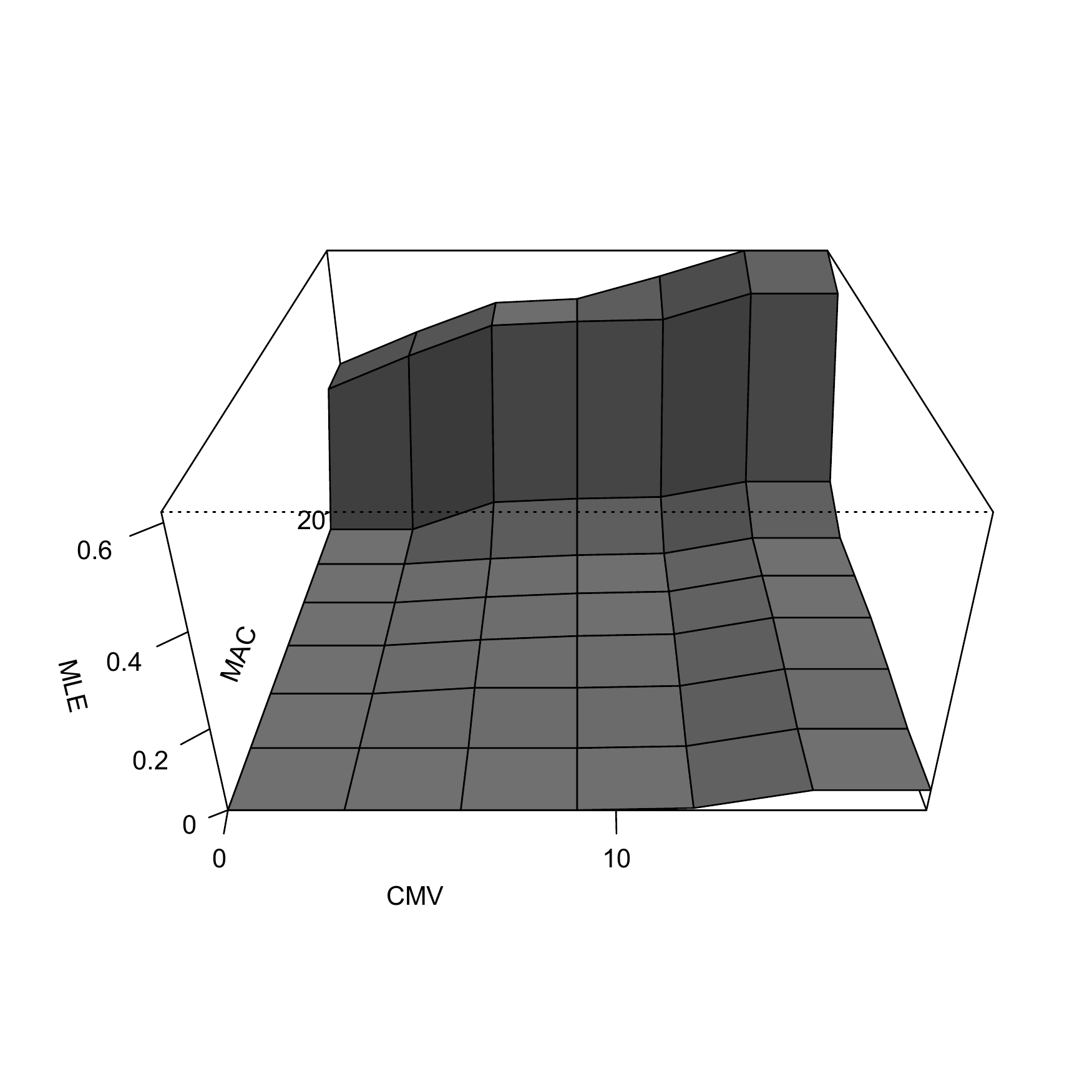}
\caption{MLE}
\label{fig:MLE_BF}
\end{subfigure}
\begin{subfigure}[b]{0.45\textwidth}
\centering
\includegraphics[width=\textwidth]{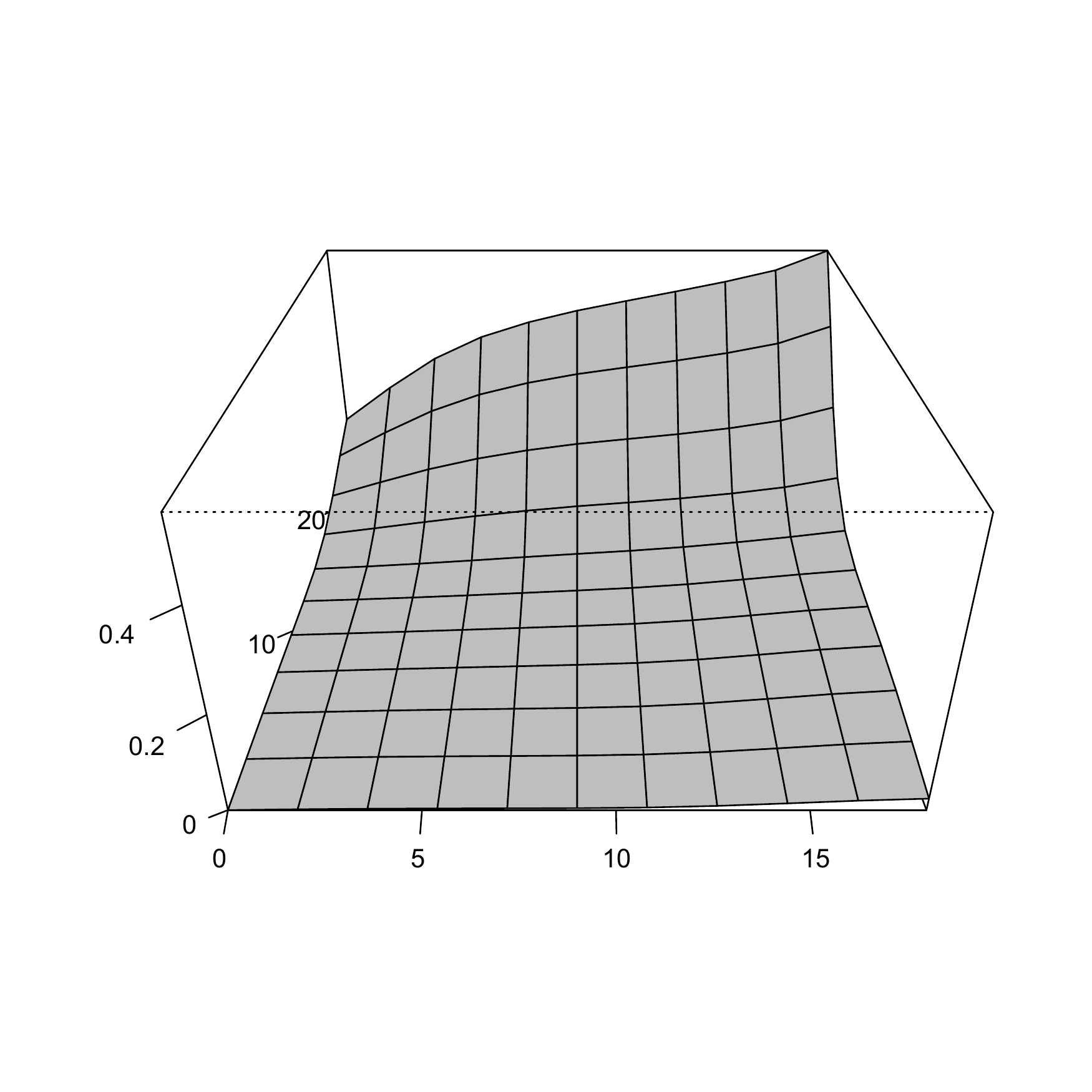}
\caption{SMLE}
\label{fig:SMLE_BF}
\end{subfigure}
\caption{MLE and SMLE for the Betensky-Finkelstein data, restricted to the interval $[0,18]\times[0,24]$.}
\label{fig:MLE_IC+MLE_IC}
\end{figure}

\begin{figure}[!ht]
\begin{subfigure}[b]{0.45\textwidth}
\centering
\includegraphics[width=\textwidth]{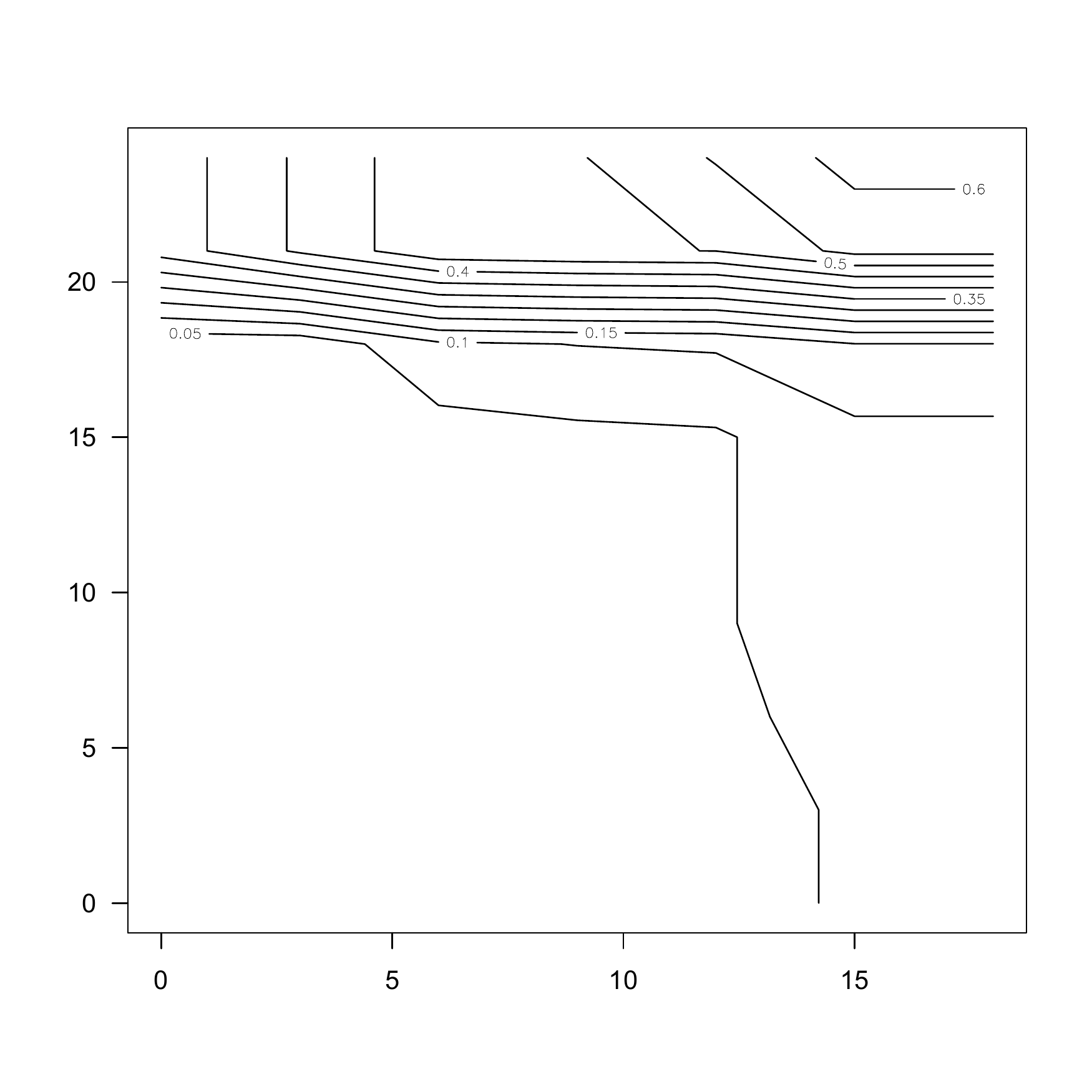}
\caption{Level plot of the MLE}
\label{fig:levels_MLE_BF}
\end{subfigure}
\begin{subfigure}[b]{0.45\textwidth}
\centering
\includegraphics[width=\textwidth]{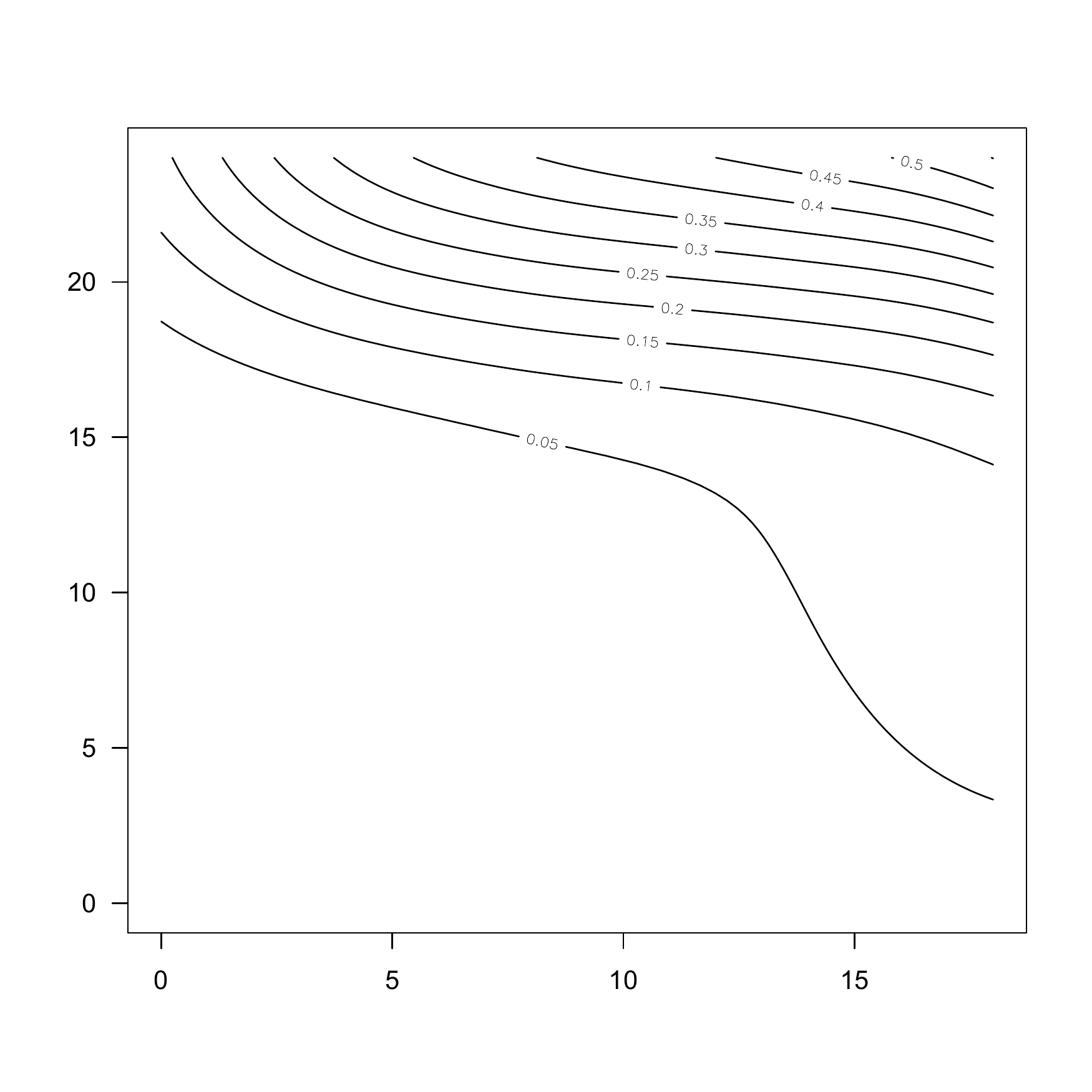}
\caption{Level plot of the SMLE}
\label{fig:levels_SMLE_BF}
\end{subfigure}
\caption{Contourplot of the MLE and SMLE for the Betensky-Finkelstein data, restricted to the interval $[0,18]\times[0,24]$.}
\label{fig:levels_MLE_IC+MLE_IC}
\end{figure}

The SMLE for bivariate interval censoring is again defined by (\ref{SMLE}). A picture of the MLE and the SMLE is shown in Figure \ref{fig:MLE_IC+MLE_IC} and the picture of the level curves in Figure \ref{fig:levels_MLE_IC+MLE_IC}. For the meaning of the codings CMV (cytomegalovirus) and MAC (mycobacterium avium complex), see \cite{Betensky:99} or \cite{sun:06}, section 7.3. Both the MSLE and the MSLE indicate that CMV shedding occurs prior to MAC colonization.

\begin{figure}[!ht]
\begin{center}
\includegraphics[scale=0.4]{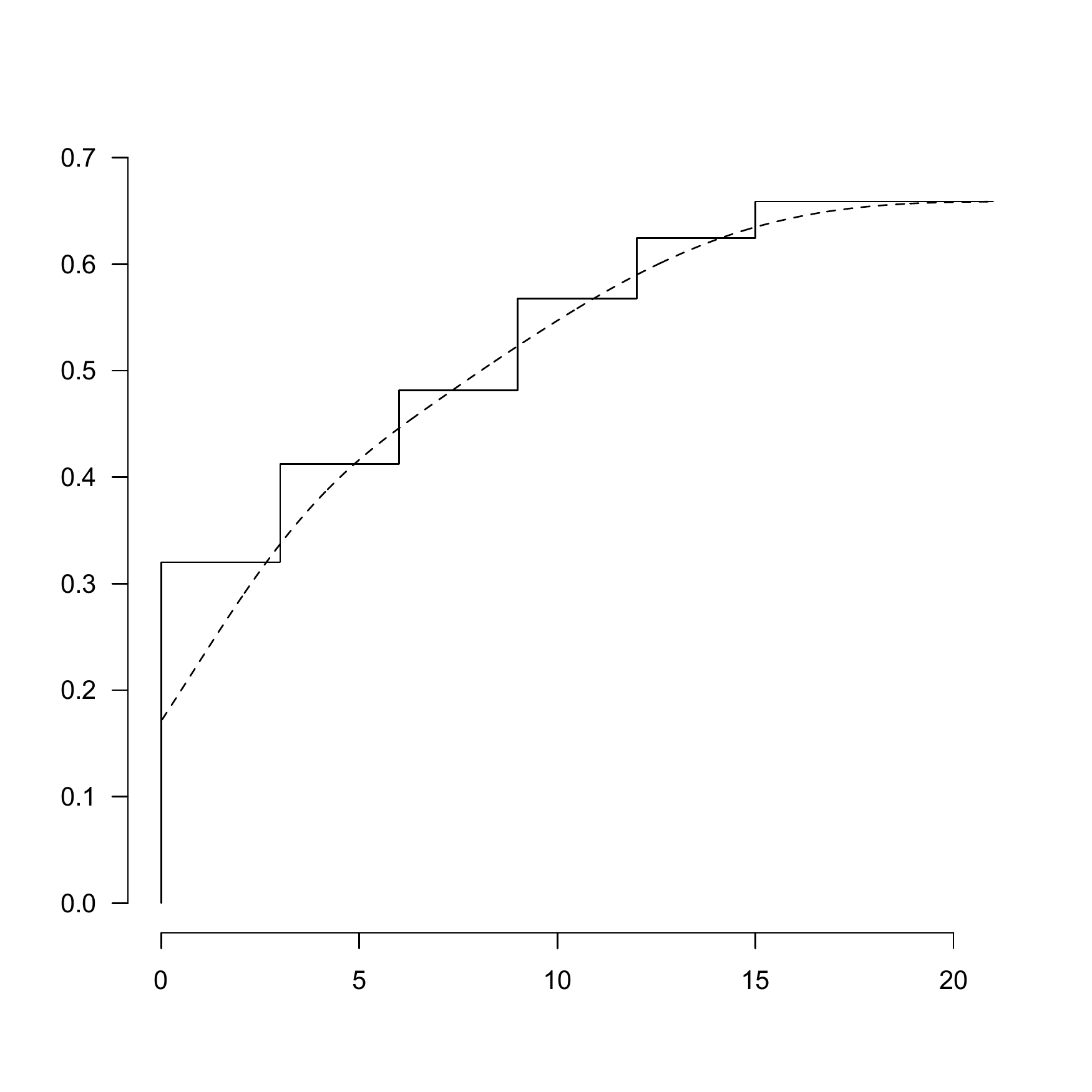}
\end{center}
\caption{The first marginal df of the data set in \cite{Betensky:99}, computed on the interval $[0,21)$ (the largest observation point on the first coordinate is $21$). The solid curve gives the first marginal df of the MLE and the dashed curve the first marginal of the SMLE, taking bandwidth $n^{-1/6}$. }
\label{fig:MLE1+SMLE1_Fink}
\end{figure}

It can be seen from this picture that the steep increase of the first marginal df of the MLE and SMLE, shown in Figure \ref{fig:MLE1+SMLE1_Fink}, is due to the `ridge' for the larger values of the second coordinate. The levels of both estimates are shown in Figure \ref{fig:levels_MLE_IC+MLE_IC}. It seems to me that the SMLE might be the more sensible estimate, also in view of the representational non-uniqueness of the MLE, which is somewhat `washed out' by the SMLE.

\section{Concluding remarks}
\label{section:conclusion}
In the preceding, three estimators for the bivariate current status model were studied: the maximum likelihood estimator (MLE), which in the simulation study was restricted to the MLE on a sieve, the smoothed maximum likelihood estimator (SMLE), obtained by integrating a kernel w.r.t.\ the masses of the MLE, and a purely discrete plug-in estimator. The SMLE and plug-in seem to attain the $n^{1/3}$ rate, with asymptotically normal distributions.

It might be somewhat surprising that the SMLE have the same rate, whereas the natural rates in the one-dimensional case are $n^{1/3}$ and $n^{2/5}$, respectively, see \cite{piet_geurt_birgit:10}. But in the bivariate case the variance is of order $1/(nh^2)$ and the bias of order $h^2$, if a bandwidth $h$ is taken in both directions and the kernel is of the usual symmetric and positive type. This makes the optimal choice of bandwidth of order $n^{-1/6}$, leading to a rate of order $n^{1/3}$.

We concentrated on compact support, but since we focused on local estimation, this restriction does not seem essential, except for the SMLE, since there the boundary played  an important role. Further research on this matter seems needed.

It is also possible to attain again the local $n^{2/5}$ rate  in the bivariate case, but then one has to take recourse to higher order kernels $K$ with the property
$$
\int u^2 K(u)\,du=0.
$$
One also has to take bandwidths of order $n^{-1/10}$ instead of order $n^{-1/6}$ in this case, which makes the judicious choice of boundary kernels even more important. Moreover, one has to strengthen the conditions of the theorems to the existence of $4$th derivatives. However, if one is willing to do that, it is easy to strengthen Theorem \ref{th:SMLE} by letting the kernel $\IK$ be based on, for example, the kernel
$$
K_1(u)=\tfrac{315}{512}\left(1-u^2\right)^3(11u^2-3)1_{[-1,1]}(u)
$$
which is the fourth order kernel corresponding to the Triweight kernel
$$
K(u)=\tfrac{35}{32}\left(1-u^2\right)^31_{[-1,1]}(u).
$$
But one loses the property that the resulting estimate is necessarily a distribution function, since the kernel $K_1$ is no longer positive. For the choice of higher order kernels, see, e.g., \cite{Jones:97}.

\section{Appendix}
\label{section:appendix}

{\bf Proof of Theorem \ref{th:plug-in}}. We use the representation
\begin{align}
\label{decomposition}
&\frac{\int_{A_n} \d_1\d_2\,d\P_n(v,w,\d_1,\d_2)}{\int_{A_n} d\G_n(v,w)}-F_0(t,u)\nonumber\\
&=\frac{\int_{A_n} \d_1\d_2\,d\P_n(v,w,\d_1,\d_2)-E\left\{\int_{A_n} \d_1\d_2\,d\P_n(v,w,\d_1,\d_2)\bigm|(T_i,U_i),\,i=1,\dots,n\right\}}{\int_{A_n} d\G_n(v,w)}\nonumber\\
&\qquad\qquad+\frac{E\left\{\int_{A_n} \d_1\d_2\,d\P_n(v,w,\d_1,\d_2)\bigm|(T_i,U_i),\,i=1,\dots,n\right\}}{\int_{A_n} d\G_n(v,w)}-\frac{n^{-1}\sum_{i=1}^n F_0(t,u)1_{A_n}(T_i,U_i)}{\int_{A_n} d\G_n(v,w)}\,.
\end{align}
The numerator of the first term on the right-hand side can be written:
\begin{align*}
n^{-1}\sum_{i=1}^n\left\{\dd_{i1}\dd_{i2}-F_0(T_i,U_i)\right\}1_{A_n}(T_i,U_i).
\end{align*}
This is the sum of i.i.d.\ random variables, and
$$
\mbox{var}\left(\left\{\dd_{11}\dd_{12}-F_0(T_1,U_1)\right\}1_{A_n}(T_1,U_1)\right)
\sim 4n^{-1/3}g(t,u)F_0(t,u)\left\{1-F_0(t,u)\right\},\,n\to\infty.
$$
Hence:
\begin{align*}
n^{1/3}\frac{\int_{A_n} \d_1\d_2\,d\P_n(v,w,\d_1,\d_2)-E\left\{\int_{A_n} \d_1\d_2\,d\P_n(v,w,\d_1,\d_2)\bigm|(T_i,U_i),\,i=1,\dots,n\right\}}{\int_{A_n} d\G_n(v,w)}
\stackrel{{\cal D}}\longrightarrow N(0,\s^2),
\end{align*}
where
$$
\s^2=\frac{F_0(t,u)\left\{1-F_0(t,u)\right\}}{4g(t,u)}\,.
$$

The numerator of the second term on the right-hand side of (\ref{decomposition}) can be written:
\begin{align*}
n^{-1}\sum_{i=1}^n\left\{F_0(T_i,U_i)-F_0(t,u)\right\}1_{A_n}(T_i,U_i).
\end{align*}
Note that, putting $h=h_n\sim n^{-1/6}$,
\begin{align*}
&n^{-1}\sum_{i=1}^n E\left\{F_0(T_i,U_i)-F_0(t,u)\right\}1_{A_n}(T_i,U_i)
=\int_{\max\{|v-t|,|w-u|\}\le h}\left\{F_0(v,w)-F_0(t,u)\right\}\,g(v,w)\,dv\,dw\\
&=\partial_1F_0(t,u)\int_{\max\{|v-t|,|w-u|\}\le h}(v-t)\,g(v,w)\,dv\,dw\\
&\qquad+\partial_2F_0(t,u)\int_{\max\{|v-t|,|w-u|\}\le h}(w-u)\,g(v,w)\,dv\,dw\\
&\qquad+\tfrac12\partial_1^2F_0(t,u)\int_{\max\{|v-t|,|w-u|\}\le h}(v-t)^2\,g(v,w)\,dv\,dw\\
&\qquad+\tfrac12\partial_2^2F_0(t,u)\int_{\max\{|v-t|,|w-u|\}\le h}(w-u)^2\,g(v,w)\,dv\,dw
+o\left(n^{-2/3}\right)\\
&=\partial_1F_0(t,u)\partial_1g(t,u)\int_{\max\{|v-t|,|w-u|\}\le h}(v-t)^2\,dv\,dw\\
&\qquad+\partial_2F_0(t,u)\partial_2g(t,u)\int_{\max\{|v-t|,|w-u|\}\le h}(w-u)^2\,dv\,dw\\
&\qquad+\tfrac12\partial_1^2F_0(t,u)\int_{\max\{|v-t|,|w-u|\}\le h}(v-t)^2\,g(v,u)\,dv\,dw\\
&\qquad+\tfrac12\partial_2^2F_0(t,u)\int_{\max\{|v-t|,|w-u|\}\le h}(w-u)^2\,g(v,u)\,dv\,dw
+o\left(n^{-2/3}\right)\\
&=\tfrac43\left\{\partial_1F_0(t,u)\partial_1g(t,u)+\partial_2F_0(t,u)\partial_2g(t,u)\right\}h^4\\
&\qquad\qquad+\tfrac23\left\{\partial_1^2F_0(t,u)+\partial_2^2F_0(t,u)\right\}g(t,u)h^4+o\left(n^{-2/3}\right).
\end{align*}
Moreover,
\begin{align*}
&\mbox{var}\left(n^{-1}\sum_{i=1}^n\left\{F_0(T_i,U_i)-F_0(t,u)\right\}1_{A_n}(T_i,U_i)\right)\\
&=O\left(n^{-1}\int_{\max\{|t-t_0|,|u-u_0|\}\le h}\left\{F_0(t,u)-F_0(t_0,u_0)\right\}^2\,g(t,u)\,dt\,du\right)\\
&=O\left(n^{-5/6}\right).
\end{align*}
The result now follows.\eop

\bibliographystyle{amsplain}
\bibliography{cupbook}

\end{document}